\documentclass{kybernetika}
\usepackage{booktabs}
\usepackage{flafter}
\usepackage{mathtools}
\usepackage[numbers,sort&compress]{natbib}
\hypersetup{hidelinks}

\newtheorem{theorem}{Theorem}[section]
\newtheorem{lemma}[theorem]{Lemma}
\newtheorem{proposition}[theorem]{Proposition}

\newtheorem{remark}[theorem]{Remark}
\newtheorem{example}[theorem]{Example}

\newcommand{\R}{\mathbb{R}}
\newcommand{\E}{\mathbb{E}}
\newcommand{\Prob}{\mathbb{P}}
\newcommand{\1}{\mathbf{1}}
\newcommand{\Ga}{\mathrm{Ga}}
\DeclareMathOperator{\sgn}{sgn}

\begin{document}
\pagestyle{myheadings}

\title{Worst-case bias of the independence\\ assumption for functions\\ of two ordinal variables}

\author{Konstantin V.\ Loginov}

\contact{Konstantin V.}{Loginov}{Independent researcher, Valencia, Spain}{kv.loginov13@gmail.com}

\markboth{K.\,V.\ Loginov}{Worst-case bias of the independence assumption}

\hidedoi
\maketitle

\begin{abstract}
Let ordinal grades $X$ and $Y$ be obtained by cutting a pair of random variables with copula $C$ at arbitrary thresholds, and let $g$ be a function of the grades. We study the largest error $S_C(g)$, over all thresholds, that is, over all marginal laws, made when the mean of $g(X,Y)$ is computed as if the grades were independent. This is the counterpart of the classical dependence-uncertainty problem: the copula is fixed and the discretisation is free. The error is caused entirely by the interaction of the inputs, measured by an index $\rho(g)$, the largest absolute $g$-volume of a rectangle of grades, which vanishes if and only if $g$ is additive. For every copula, $S_C(g)$ lies between $\rho(g)$ times the largest deviation of $C$ from the independence copula and $\rho(g)$ times half the total variation of their difference, and the factor one half cannot be improved. For the Gaussian copula with correlation $\lambda$ we prove that $S_C(g)$ equals $|\lambda|\rho(g)/(2\pi)$ to first order, with a sharp constant and an explicit remainder. The same first-order law holds for copula families whose first-order departure from independence is a product of unimodal factors, including the Farlie--Gumbel--Morgenstern, Ali--Mikhail--Haq, Frank, Clayton and Plackett families, and we characterise the Archimedean families with this property. An application to the TNM staging of colon cancer, with data from a clinical trial, compares the bounds with the actual error of the independence computation.
\end{abstract}

\keywords{copula, Hoeffding identity, weak dependence, ordinal data, total variation, absolute regularity, decision table}

\classification{62H05, 60E15, 62H20}

\section{INTRODUCTION}\label{sec:intro}

Let $(U,V)$ be a pair of random variables with copula $C$ \citep{Nelsen2006,Joe2014,DuranteSempi2016}, and let ordinal grades $X\in[r]=\{1,\dots,r\}$ and $Y\in[c]$ be obtained by cutting $U$ and $V$ at thresholds $u_1\le\dots\le u_{r-1}$ and $v_1\le\dots\le v_{c-1}$ in $[0,1]$: $X\le i\iff U\le u_i$ and $Y\le j\iff V\le v_j$. Then
\[
u_i=\Prob(X\le i),\qquad v_j=\Prob(Y\le j),\qquad \Prob(X\le i,Y\le j)=C(u_i,v_j),
\]
and every pair of marginal laws on $[r]$ and $[c]$ arises in this way, from the thresholds $u_i$, $v_j$ equal to the values of the two distribution functions. For discrete variables the copula of $(X,Y)$ is not unique: only its values on the grid of marginal distribution function values matter \citep{GenestNeslehova2007}; the consequences for modelling and inference are discussed by \citet{Faugeras2017} and \citet{Geenens2020}. We therefore treat $C$ as a chosen parametric model of the dependence between the latent variables, from which the joint law of the grades is built through these grid values. A function $g:[r]\times[c]\to\R$ of the grades is often evaluated from the two marginal laws as if $X$ and $Y$ were independent. We ask how large the resulting error in $\E g(X,Y)$ can be, uniformly over all thresholds, that is, over all marginal laws.

This is the counterpart of the classical problem of dependence uncertainty. There the marginals are fixed and the copula varies over a class, as in the Fr\'echet--Hoeffding bounds and their refinements \citep{PuccettiWang2015}; Hoeffding-type covariance bounds \citep{Cuadras2002,Beare2009} also take the marginals as given. Here the dependence structure is fixed and the discretisation is free. A bound uniform in the thresholds is what is needed when a scale or a table is designed before data are collected, or when one table is applied to populations with different grade distributions; when both the marginals and the copula are known, the error is simply computed from \eqref{eq:identity} below. In this reading the two-point marginals that attain our lower bounds are simply the coarsest discretisations, which collapse each grade scale to two levels. A motivating case is evaluation procedures that merge two grades through a table: decision tables in the qualitative multi-criteria method DEX \citep{BohanecRajkovic1990,Bohanec2013} and the hierarchical matrix mechanisms of complex evaluation, in which the grades of two criteria are merged by a monotone matrix and the merged grade is passed to the next level \citep{Alekseev2020,Alekseev2025}. When the grades are random, for instance because each object is assessed in a series of stochastic trials, the distribution of the merged grade is routinely computed from the marginal distributions under the assumption that the criteria are independent \citep[Definition~3]{Alekseev2025}, as in the probabilistic extension of DEX \citep{TrdinBohanec2018}. The same structure appears in medicine, where the TNM stage of a tumour is read from a table of two ordinal grades (Example~\ref{ex:colon}). If both criteria are measured in the same trials, a common trial effect makes them positively dependent, and the independence computation is biased.

Latent copula models of this kind are widely used to simulate and fit ordinal data, including the Frank copula as a latent model with an analysis of correlation attenuation \citep{Barbiero2021,PanJoe2024,Barbiero2026}. In the Gaussian case this is the latent-variable model $X=i\iff a_{i-1}<Z_1\le a_i$, $Y=j\iff b_{j-1}<Z_2\le b_j$ with $\mathrm{corr}(Z_1,Z_2)=\lambda$. Here $\lambda$ is the latent (polychoric-type) correlation, not the correlation of the observed grades, which is smaller in absolute value: for two grades split at their medians the observed correlation is $(2/\pi)\arcsin\lambda$, e.g.\ $0.19$ for $\lambda=0.3$, whereas the maximal correlation of $Z_1,Z_2$ over all functions is $|\lambda|$ \citep{Gebelein1941}. The constant $\arcsin|\lambda|/(2\pi)$ of our lower bound is the covariance of two median-split indicators, that is, $\frac14\cdot\frac2\pi\arcsin|\lambda|$; the restriction to threshold events is what replaces the maximal correlation $|\lambda|$ by $(2/\pi)\arcsin|\lambda|$. The bias
\[
B_C(g;X,Y)=\E_C\,g(X,Y)-\E_\Pi\,g(X,Y),\qquad \Pi(u,v)=uv,
\]
is, by a discrete form of Hoeffding's covariance identity (Lemma~\ref{lem:hoeffding}),
\begin{equation}\label{eq:identity}
B_C(g;X,Y)=\sum_{i=1}^{r-1}\sum_{j=1}^{c-1}\Delta_{ij}(g)\,\bigl[C(u_i,v_j)-u_iv_j\bigr],
\end{equation}
where $\Delta_{ij}(g)=g(i+1,j+1)-g(i+1,j)-g(i,j+1)+g(i,j)$ are the mixed differences of $g$. Identities of this type go back to \citet{Hoeffding1940}; for functions of bounded variation see \citet{Cuadras2002} and \citet{Beare2009}, who also derives upper bounds on the covariance through the Vitali variation $V(g)=\sum_{i,j}|\Delta_{ij}(g)|$. Functional generalisations of Hoeffding's identity are given by \citet{Lo2017}.

We are interested in the worst case
\[
S_C(g)=\sup_{X,Y}\lvert B_C(g;X,Y)\rvert ,
\]
the supremum being over all thresholds, that is, over all pairs of marginal laws on $[r]$ and $[c]$. Upper bounds of Vitali type are not sharp when $\Delta(g)$ has entries of both signs, which is the typical case for evaluation tables: an ``and''-like region of the table (where $\Delta>0$) may sit next to an ``or''-like region (where $\Delta<0$). Our main result identifies the exact first-order behaviour of $S_C(g)$. Define
\begin{equation}\label{eq:rho}
\rho(g)=\max_{I,J}\lvert \sum_{i\in I}\sum_{j\in J}\Delta_{ij}(g)\rvert ,
\end{equation}
where $I\subseteq[r-1]$ and $J\subseteq[c-1]$ range over nonempty discrete intervals $\{i_1,i_1+1,\dots,i_2\}$. Intervals arise because, for increasing thresholds, the superlevel sets of unimodal weights are intervals of grades (Lemma~\ref{lem:layer}). Thus $\rho(g)$ is the largest absolute sum of a contiguous submatrix of $\Delta(g)$ (a ``maximum subarray'' in the sense of \citealp{Bentley1984}), an interval analogue of the cut norm of \citet{FriezeKannan1999}. Since a rectangle sum of mixed differences telescopes,
\[
\rho(g)=\max_{i<i',\,j<j'}\lvert g(i',j')-g(i',j)-g(i,j')+g(i,j)\rvert ,
\]
the largest absolute $g$-volume of a rectangle of grades, that is, the largest interaction contrast of $g$; in particular $\rho(g)=0$ if and only if $g$ is additive, $g(i,j)=a_i+b_j$. For tables with a few grades this direct form is the simplest way to compute $\rho$. Clearly $|N(g)|\le\rho(g)\le\|\Delta(g)\|_\square\le V(g)$, where $N(g)=\sum_{i,j}\Delta_{ij}(g)=g(r,c)-g(r,1)-g(1,c)+g(1,1)$ is the corner contrast and $\|\cdot\|_\square$ is the cut norm, the same maximum over arbitrary index sets (Section~\ref{sec:discussion}). The index is homogeneous of degree one in $g$ and does not change when a constant is added to $g$. Since every rectangle sum of mixed differences is a signed sum of four values of $g$ with two signs of each kind, $\rho(g)\le2(\max g-\min g)$. Non-constant tables with different output ranges can therefore be compared through the scale-free ratio $\rho(g)/(\max g-\min g)\in[0,2]$; for instance it equals $1/3$ for the matrix $M$ of Example~\ref{ex:matrix} and $1$ for $\min(i,j)$ on four grades. Each value $g(k,l)$ enters a rectangle sum of mixed differences at most once, with coefficient $\pm1$, so changing one cell of the table by $\delta$ changes $\rho$ by at most $|\delta|$; this bound is attained.

\textbf{Results.} Since $\rho(g)=0$ if and only if $g$ is additive, the error of the independence computation is caused entirely by the interaction of the two inputs, and our results show that, under weak dependence, its worst case is governed by the largest interaction contrast $\rho(g)$. For every copula, $\rho(g)\varepsilon(C)\le S_C(g)\le V(g)\varepsilon(C)$ with $\varepsilon(C)=\sup_{u,v}|C(u,v)-uv|$ (Proposition~\ref{prop:general}). At every strength of dependence the variation $V(g)$ can be replaced by $\rho(g)$ at the price of a different dependence measure: $S_C(g)\le\rho(g)\tau(C)/2$, where $\tau(C)$ is the total variation norm of the difference between $C$ and the independence copula, that is, twice the total variation distance and twice the absolute regularity coefficient of the latent pair (Theorem~\ref{thm:tv}). The constant $1/2$ is sharp, and the best constant $K(C)=\sup_gS_C(g)/\rho(g)$ always lies between $\tau(C)/4$ and $\tau(C)/2$. For the Gaussian copula $C^{\Ga}_\lambda$ we show (Theorem~\ref{thm:gauss}) that
\[
\frac{\arcsin|\lambda|}{2\pi}\,\rho(g)\;\le\;S_\lambda(g)\;\le\;\frac{|\lambda|}{2\pi}\,\rho(g)+V(g)\,R(|\lambda|),\qquad |\lambda|<1,
\]
with an explicit remainder $R(\lambda)=\lambda^2/(4\pi\mathrm e)+O(\lambda^3)$ that is built on a pointwise sharp bound and is asymptotically sharp in the sense that the supremum over thresholds of the per-cell remainder is $R(\lambda)(1+o(1))$ as $\lambda\to0$. A second upper bound, $S_\lambda(g)\le\rho(g)G(|\lambda|)$ with an explicit $G$, follows from the tetrachoric series and does not involve $V(g)$. Consequently $S_\lambda(g)=|\lambda|\rho(g)/(2\pi)+O(\lambda^2)$: to first order, the worst-case sensitivity of a table to weak Gaussian dependence is a property of the table alone, equal to $\rho(g)/(2\pi)$ per unit of correlation. The lower bound is attained by two-point marginals that split each variable at the median on the border of the optimal rectangle. Theorem~\ref{thm:families} extends the first-order law to copula families $C_\theta$ with $C_\theta(u,v)=uv+(\theta-\theta_0)A(u)A(v)+O((\theta-\theta_0)^2)$ and unimodal $A\ge0$; then $S_{C_\theta}(g)=|\theta-\theta_0|\kappa\,\rho(g)+O((\theta-\theta_0)^2)$ with $\kappa=(\max A)^2$. The ranking of tables by sensitivity to weak dependence is therefore the same for all such families. Relative to the existing covariance inequalities the new elements are threefold. First, Hoeffding-type bounds \citep{Cuadras2002,Beare2009} bound the bias for given marginals by the variation $V(g)$ times a dependence measure; summed over cells, this is the $V(g)R(|\lambda|)$ part of our upper bound. We optimise over all marginals, replace $V(g)$ by $\rho(g)$ in the leading term through the unimodality of the weights (Lemma~\ref{lem:layer}), and obtain a matching lower bound. Second, the exact first-order constant is governed by the new index $\rho(g)\le V(g)$, and the matching lower bound is attained constructively by two-point marginals. Third, the non-asymptotic bound of Theorem~\ref{thm:tv} shows that $\rho(g)$ governs the worst case for every copula and every strength of dependence, up to a factor depending on the copula only. We are not aware of earlier results that identify this worst case over the marginals.

\textbf{Relation to Hoeffding's extremal principle.} \citet{Hoeffding1955} showed that extrema of $\E K(X_1,\dots,X_n)$ over independent $X_j$ subject to moment constraints are attained by distributions with few support points. Our problem has a similar flavour (an extremum over marginals, attained by two-point laws) but concerns the difference between a dependent and an independent model. The reduction to two-point laws is exact for the lower bound and for the first-order term, but not at finite dependence strength (Section~\ref{sec:numerics}).

\section{THE BIAS IDENTITY AND TWO ELEMENTARY LEMMAS}\label{sec:prelim}

\begin{lemma}[Discrete Hoeffding identity]\label{lem:hoeffding}
For every $g:[r]\times[c]\to\R$ and every joint law of $(X,Y)$ with marginal distribution functions $u_i,v_j$ and copula $C$, identity \eqref{eq:identity} holds. Equivalently, $B_C(g)=\sum_{i,j}\Delta_{ij}(g)\,[\Prob(X>i,Y>j)-\Prob(X>i)\Prob(Y>j)]$.
\end{lemma}

\begin{Proof}
For $(x,y)\in[r]\times[c]$ the double sum of $\Delta_{ij}(g)$ over $i<x$, $j<y$ telescopes, so
$g(x,y)=g(x,1)+g(1,y)-g(1,1)+\sum_{i,j}\Delta_{ij}(g)\,\1\{x>i\}\1\{y>j\}$.
Take expectations under the dependent and the independent model and subtract. The first three terms depend on the marginals only and cancel. Finally $\Prob(X>i,Y>j)=1-u_i-v_j+C(u_i,v_j)$ and $\Prob(X>i)\Prob(Y>j)=1-u_i-v_j+u_iv_j$.
\end{Proof}

The quantities $C(u_i,v_j)-u_iv_j$ are the orthant covariances of the indicators $\1\{X>i\}$, $\1\{Y>j\}$. If $u_i\in\{0,1\}$ or $v_j\in\{0,1\}$ they vanish.

Call a finite sequence $w=(w_1,\dots,w_n)$ \emph{unimodal} if $w_k\ge\min(w_i,w_l)$ whenever $i<k<l$, that is, if it is nondecreasing up to some index and nonincreasing afterwards.

For $w\in\R^m$ put
\[
\nu(w)=\frac12\sum_{k=0}^{m}|w_{k+1}-w_k|,\qquad w_0=w_{m+1}=0,
\]
half the total variation of $w$ extended by zeros at both ends.

\begin{lemma}[Layer-cake bound]\label{lem:layer}
Let $D\in\R^{m\times n}$, $w\in\R^m$ and $z\in\R^n$, and let $\rho(D)$ be the largest absolute sum of a contiguous submatrix of $D$, as in \eqref{eq:rho}. Then $|w^\top Dz|\le\rho(D)\,\nu(w)\,\nu(z)$. If $w\in[0,1]^m$ is unimodal, then $\nu(w)=\max_iw_i\le1$; in particular $|w^\top Dz|\le\rho(D)$ for unimodal $w\in[0,1]^m$, $z\in[0,1]^n$. The inequality is sharp: equality holds when $w$ and $z$ are indicators of the intervals attaining $\rho(D)$.
\end{lemma}

\begin{Proof}
First let $w\ge0$, and let $0=\ell_0<\ell_1<\dots<\ell_K$ be its distinct positive values together with $0$. Then $w=\sum_{t=1}^K(\ell_t-\ell_{t-1})\1_{\{w\ge\ell_t\}}$. If the superlevel set $\{w\ge\ell_t\}$ consists of $k_t$ maximal sets of consecutive integers, then each of them contributes one up-step and one down-step through the level $\ell_t$ to the sequence $(0,w_1,\dots,w_m,0)$, so $\sum_t(\ell_t-\ell_{t-1})k_t=\nu(w)$. Hence $w$ is a combination of indicators of discrete intervals whose coefficients have total absolute value $\nu(w)$. For general $w$ apply this to $w^\pm=\max(\pm w,0)$; since $w^+_k$ and $w^-_k$ are never both positive, $|w_{k+1}-w_k|=|w^+_{k+1}-w^+_k|+|w^-_{k+1}-w^-_k|$, so $\nu(w)=\nu(w^+)+\nu(w^-)$. Decomposing $z$ in the same way, $|w^\top Dz|$ is at most $\nu(w)\nu(z)$ times the largest $|\1_I^\top D\1_J|$ over intervals $I,J$, which is $\rho(D)$. If $w\in[0,1]^m$ is unimodal, each superlevel set is a single interval (if $i<k<l$ and $w_i,w_l\ge\ell_t$, then $w_k\ge\min(w_i,w_l)\ge\ell_t$), so $k_t=1$ and $\nu(w)=\ell_K=\max_iw_i$.
\end{Proof}

Lemma~\ref{lem:layer} is a discrete analogue of the Koksma--Hlawka inequality \citep{KuipersNiederreiter1974}: the integral of a function of bounded variation against a signed measure is bounded by the variation times the discrepancy, here with $\nu(w)\nu(z)$ in the role of the variation of $w\otimes z$ and $\rho(D)$, the largest mass of a box, in that of the discrepancy.

\textbf{A general frame.} For a copula $C$ let $\varepsilon(C)=\sup_{u,v}|C(u,v)-uv|$; this is a quarter of the $L_\infty$ measure of dependence of \citet{SchweizerWolff1981}. Let $P(g)$ and $Q(g)$ be the sums of the positive parts and of the negative parts of the entries of $\Delta(g)$.

\begin{proposition}\label{prop:general}
For every copula $C$ and every $g$,
\[
\rho(g)\,\varepsilon(C)\;\le\;S_C(g)\;\le\;V(g)\,\varepsilon(C),
\]
and the lower bound is attained by two-point marginals; in particular, $S_C(g)=|N(g)|\,\varepsilon(C)$ if $\Delta(g)$ has constant sign, since then $\rho(g)=V(g)=|N(g)|$. If $C\ge\Pi$ or $C\le\Pi$ pointwise (positive or negative quadrant dependence), then $S_C(g)$ is at most $\max\{P(g),Q(g)\}\,\varepsilon(C)$.
\end{proposition}

\begin{Proof}
Write $c_{ij}=C(u_i,v_j)-u_iv_j$. By \eqref{eq:identity}, $|B_C|\le\sum|\Delta_{ij}||c_{ij}|\le V(g)\varepsilon(C)$; if all $c_{ij}$ have one sign, $B_C$ lies between $-Q\varepsilon$ and $P\varepsilon$ or between $-P\varepsilon$ and $Q\varepsilon$. For the lower bound, let $I=\{i_1,\dots,i_2\}$ and $J=\{j_1,\dots,j_2\}$ attain \eqref{eq:rho}, and let $(u^\ast,v^\ast)$ attain $\varepsilon(C)$ (the supremum of a continuous function on $[0,1]^2$; if $C\ne\Pi$, then $u^\ast,v^\ast\in(0,1)$ because $C-\Pi$ vanishes on the boundary). Let $X=i_1$ with probability $u^\ast$ and $X=i_2+1$ otherwise, and $Y=j_1$ with probability $v^\ast$ and $Y=j_2+1$ otherwise. Then $u_i=u^\ast$ for $i\in I$, $v_j=v^\ast$ for $j\in J$, and all other $u_i,v_j$ lie in $\{0,1\}$, where $c_{ij}=0$. Hence $B_C=\bigl(\sum_{I\times J}\Delta_{ij}\bigr)\bigl(C(u^\ast,v^\ast)-u^\ast v^\ast\bigr)$ and $|B_C|=\rho(g)\varepsilon(C)$. For constant-sign $\Delta$, $\rho=V=|N|=\max(P,Q)$.
\end{Proof}

For constant-sign $\Delta$, that is, for supermodular or submodular $g$, the sign of the bias under positive quadrant dependence is classical: $\E_Cg\ge\E_\Pi g$ for supermodular $g$ and $C\ge\Pi$ \citep{Lehmann1966,CambanisSimonsStout1976,Tchen1980}. Proposition~\ref{prop:general} reduces the problem to the question when $V(g)$ can be replaced by $\rho(g)$ in the upper bound. Theorems~\ref{thm:gauss} and~\ref{thm:families} show that this is possible to first order in the dependence parameter, and the Gumbel example of Section~\ref{sec:families} shows that in general $\rho(g)$ determines the first-order term only up to a constant factor.

\section{A BOUND AT EVERY STRENGTH OF DEPENDENCE}\label{sec:tv}

Proposition~\ref{prop:general} brackets $S_C(g)$ between $\rho(g)\varepsilon(C)$ and $V(g)\varepsilon(C)$, and the two indices can differ by a large factor. In this section we show that, at the price of a different dependence measure, $\rho(g)$ also controls the upper bound, for every copula and every strength of dependence.

For a copula $C$ let $\mu_C$ be the probability measure on $[0,1]^2$ with distribution function $C$, and let
\[
\tau(C)=|\mu_C-\mu_\Pi|\bigl([0,1]^2\bigr)
\]
be the total variation of the signed measure $\mu_C-\mu_\Pi$, that is, the total variation norm, which is twice the total variation distance $\sup_A|\mu_C(A)-\mu_\Pi(A)|$; if $C$ has a density $c$, then $\tau(C)=\int_0^1\!\int_0^1|c(u,v)-1|\,du\,dv$. For a pair $(U,V)$ with copula $C$, $\tau(C)/2$ is the absolute regularity ($\beta$-mixing) coefficient $\beta(\sigma(U),\sigma(V))$ \citep{Bradley2005}; for copula-based Markov chains it is the coefficient studied by \citet{LonglaPeligrad2012}. Since the four quadrants at $(u,v)$ carry the masses $\pm(C(u,v)-uv)$ under $\mu_C-\mu_\Pi$, we have $\varepsilon(C)\le\tau(C)/4$. We measure the overall sensitivity of the copula by the best constant
\[
K(C)=\sup\Bigl\{\frac{S_C(g)}{\rho(g)}:\ r,c\ge2,\ g:[r]\times[c]\to\R\ \text{not additive}\Bigr\};
\]
by Proposition~\ref{prop:general}, $K(C)\ge\varepsilon(C)$.

\begin{lemma}[Additive approximation]\label{lem:additive}
Let $g:[r]\times[c]\to\R$ and $m=\min(r,c)$. There are $a\in\R^r$ and $b\in\R^c$ with
\[
\max_{i,j}\lvert g(i,j)-a_i-b_j\rvert \le\frac{m-1}{2m}\,\rho(g).
\]
The constant $(m-1)/(2m)$ cannot be improved.
\end{lemma}

\begin{Proof}
By symmetry we may assume $r=m$. Put $b_j=\frac1m\sum_{i'}g(i',j)$. For fixed $i$,
\[
g(i,j)-b_j=\frac1m\sum_{i'\ne i}\bigl(g(i,j)-g(i',j)\bigr).
\]
For each $i'$ the function $j\mapsto g(i,j)-g(i',j)$ has oscillation $\max_j-\min_j$ at most $\rho(g)$, because its increments $g(i,j')-g(i',j')-g(i,j)+g(i',j)$ are $g$-volumes of rectangles. Hence $j\mapsto g(i,j)-b_j$ has oscillation at most $(m-1)\rho(g)/m$, and choosing $a_i$ as the midpoint of its range gives the bound. Sharpness is shown in the proof of Theorem~\ref{thm:tv}.
\end{Proof}

The quantity $\min_{a,b}\max_{i,j}|g(i,j)-a_i-b_j|$ is the classical problem of best uniform approximation by sums of functions of one variable \citep{DilibertoStraus1951}. The constant $2-2/m=2\cdot\frac{m-1}{m}$ is the projection constant of the hyperplane of zero-sum vectors in $\ell_1^m$ \citep[Lemma~3]{BlatterCheney1974}, and the proof above is the projection with the uniform weights.

\begin{theorem}\label{thm:tv}
For every copula $C$ and every $g:[r]\times[c]\to\R$, with $m=\min(r,c)$,
\[
S_C(g)\;\le\;\frac{m-1}{2m}\,\rho(g)\,\tau(C)\;\le\;\frac12\,\rho(g)\,\tau(C),
\]
with strict inequality in the second step whenever $\rho(g)\tau(C)>0$.
Moreover $\tau(C)/4\le K(C)\le\tau(C)/2$, and both constants are sharp: $K(C)=\tau(C)/4$ for every copula $C(u,v)=uv+\theta A(u)B(v)$ with absolutely continuous $A,B$, whereas $\sup_{C\ne\Pi}K(C)/\tau(C)=\frac12$.
\end{theorem}

\begin{Proof}
\emph{Upper bound.} Fix the thresholds $0=u_0\le u_1\le\dots\le u_r=1$ and $0=v_0\le\dots\le v_c=1$, and let $q_{ij}$ be the mass that $\mu_C-\mu_\Pi$ gives to the cell $(u_{i-1},u_i]\times(v_{j-1},v_j]$, so that $\Prob_C(X=i,Y=j)-\Prob_\Pi(X=i,Y=j)=q_{ij}$ and $B_C(g)=\sum_{i,j}g(i,j)q_{ij}$. The two models have the same marginals, so every row sum and every column sum of $(q_{ij})$ vanishes. Hence $B_C(g)=\sum_{i,j}\bigl(g(i,j)-a_i-b_j\bigr)q_{ij}$ for all $a,b$, and Lemma~\ref{lem:additive} gives $|B_C(g)|\le\frac{m-1}{2m}\rho(g)\sum_{i,j}|q_{ij}|\le\frac{m-1}{2m}\rho(g)\tau(C)$.

\emph{Lower bound for $K$.} For a grid of thresholds as above put $g(i,j)=\sgn q_{ij}$. Each $g$-volume of a rectangle is a signed sum of four numbers in $[-1,1]$, so $\rho(g)\le4$, while $B_C(g)=\sum|q_{ij}|$. Hence $K(C)\ge\frac14\sum_{i,j}|q_{ij}|$ for every grid. Let $P$ be a positive set of a Hahn decomposition of $\nu=\mu_C-\mu_\Pi$. Since finite unions of rectangles approximate Borel sets in $|\nu|$-measure, for every $\delta>0$ there is a finite union $F$ of rectangles with $|\nu|(P\triangle F)<\delta$; refining the grid until $F$ is a union of cells gives $\sum|q_{ij}|\ge\nu(F)-\nu(F^c)\ge\tau(C)-4\delta$. Thus $K(C)\ge\tau(C)/4$.

\emph{Product copulas.} Let $C(u,v)=uv+\theta A(u)B(v)$. Then $q_{ij}=\theta\alpha_i\beta_j$ with $\alpha_i=A(u_i)-A(u_{i-1})$ and $\beta_j=B(v_j)-B(v_{j-1})$, and $\sum_i\alpha_i=\sum_j\beta_j=0$ because $A$ and $B$ vanish at $0$ and $1$. A vector with zero sum is a nonnegative combination $\sum_ks_k(e_{i_k}-e_{i'_k})$ of differences of unit vectors with $\sum_ks_k=\frac12\|\alpha\|_1$ (pair its positive part with its negative part); write $\beta=\sum_lt_l(e_{j_l}-e_{j'_l})$ in the same way. Then
\[
B_C(g)=\theta\sum_{k,l}s_kt_l\bigl[g(i_k,j_l)-g(i_k,j'_l)-g(i'_k,j_l)+g(i'_k,j'_l)\bigr],
\]
each bracket is at most $\rho(g)$ in absolute value, and so $|B_C(g)|\le\frac14|\theta|\,\|\alpha\|_1\|\beta\|_1\rho(g)=\frac14\rho(g)\sum|q_{ij}|\le\frac14\rho(g)\tau(C)$. Hence $K(C)=\tau(C)/4$.

\emph{The constant $\frac12$.} Let $m\ge2$, index the cells of the uniform $m\times m$ grid by $i,j\in\{0,\dots,m-1\}$, and let $q_{ij}=1$ if $j=i$, $q_{ij}=-1$ if $j\equiv i+1\pmod m$, and $q_{ij}=0$ otherwise. Every row and column of $q$ contains one $+1$ and one $-1$. The checkerboard copula $C_m$ with density $1+q_{ij}\in\{0,1,2\}$ on cell $(i,j)$ therefore has uniform marginals, and $\tau(C_m)=2m/m^2=2/m$. Let $g(i,j)=\1\{1\le j\le i\}$. For $i<i'$ the difference $g(i',\cdot)-g(i,\cdot)$ takes values in $\{0,1\}$, so $\rho(g)=1$. With the uniform marginals on $m$ grades (thresholds $k/m$), $B_{C_m}(g)=m^{-2}\sum_{i,j}g(i,j)q_{ij}=m^{-2}\sum_{i=0}^{m-1}\bigl[g(i,i)-g(i,i+1\bmod m)\bigr]=(m-1)/m^2$, since $g(i,i)=1$ for $i\ge1$ and $g(i,i+1\bmod m)=0$ for all $i$. Thus $S_{C_m}(g)\ge\frac{m-1}{2m}\rho(g)\tau(C_m)$: the first bound of the theorem, and hence Lemma~\ref{lem:additive}, is sharp, and $K(C_m)/\tau(C_m)\ge\frac{m-1}{2m}\to\frac12$.
\end{Proof}

Combined with Proposition~\ref{prop:general}, Theorem~\ref{thm:tv} gives, for every copula and every table,
\[
\rho(g)\,\varepsilon(C)\;\le\;S_C(g)\;\le\;\rho(g)\,\frac{\tau(C)}{2}:
\]
the table enters the worst case only through $\rho(g)$, at every strength of dependence, up to the factor $\tau(C)/(2\varepsilon(C))\ge2$ that depends on the copula alone. For the Gaussian copula this factor equals $2.04$, $2.31$ and $2.89$ at $|\lambda|=0.2$, $0.6$ and $0.9$ (numerical quadrature of $\tau$). The tetrachoric bound of Theorem~\ref{thm:gauss} is sharper for small and moderate $|\lambda|$, with factors $1.17$ and $1.77$ at $0.2$ and $0.6$, but Theorem~\ref{thm:tv} is sharper for strong dependence ($2.89$ against $3.66$ at $|\lambda|=0.9$), and it applies to every copula, including those, such as the Gumbel and Joe families, whose first-order term is not of product form. For product copulas $uv+\theta A(u)B(v)$ with unimodal $A,B$ the bound $K(C)=\tau(C)/4=|\theta|\kappa$ agrees with Theorem~\ref{thm:families}, since then $\mathrm{TV}(A)=2\max A$. Which copulas satisfy $K(C)=\tau(C)/4$ in general is open. By the argument used for product copulas, a sufficient condition is that on every grid $(q_{ij})$ is a sum of terms $w\,(e_i-e_{i'})(e_j-e_{j'})^\top$, $w\ge0$, whose signs agree with those of $q$. For the Gaussian copula with small $|\lambda|$ one can show that $K(C)>\tau(C)/4$; this will be reported elsewhere.

\section{THE GAUSSIAN COPULA}\label{sec:gauss}

Let $\varphi,\Phi$ be the standard normal density and distribution function, and let $\varphi_2(a,b;s)$, $\Phi_2(a,b;s)$ be the standard bivariate normal density and distribution function with correlation $s\in(-1,1)$. The Gaussian copula is $C^{\Ga}_\lambda(u,v)=\Phi_2(\Phi^{-1}(u),\Phi^{-1}(v);\lambda)$. We write $S_\lambda(g)$ for $S_{C^{\Ga}_\lambda}(g)$.

\begin{lemma}[Plackett's identities]\label{lem:plackett}
For $a,b\in\R$ and $|s|<1$,
$\partial_s\Phi_2(a,b;s)=\varphi_2(a,b;s)$ and $\partial_s\varphi_2(a,b;s)=\partial_a\partial_b\varphi_2(a,b;s)$.
\end{lemma}

These are classical \citep{Plackett1954}; the first one underlies Price's theorem \citep{Price1958}. Both follow from a direct computation with the explicit density.

\begin{lemma}[Sharp bound on the mixed derivative]\label{lem:h}
For $s\in[0,1)$ let
\[
h(s)=\frac{1}{2\pi\sqrt{1-s^2}}\max\Bigl\{h_+(s),\;\frac{\mathrm e^{-(1+2s)/(1+s)}}{1-s}\Bigr\},
\]
where
\[
h_+(s)=\begin{cases}\dfrac{\mathrm e^{-(1-2s)/(1-s)}}{1+s}, & s<\tfrac12,\\[2ex] \dfrac{s}{1-s^2}, & s\ge\tfrac12,\end{cases}
\]
and $h(-s)=h(s)$. Then $\sup_{a,b\in\R}|\partial_a\partial_b\varphi_2(a,b;s)|=h(s)$. In particular $h(0)=\mathrm e^{-1}/(2\pi)$. Moreover, $h$ is continuous and nondecreasing on $[0,1)$.
\end{lemma}

\begin{Proof}
Since $\varphi_2(-a,b;-s)=\varphi_2(a,b;s)$, the mixed derivative only changes sign under $s\mapsto-s$, so we may assume $s\ge0$. Put $q=1-s^2$. Differentiation gives
$\partial_a\partial_b\varphi_2=\varphi_2\cdot\bigl[(a-sb)(b-sa)/q^2+s/q\bigr]$.
Rotate coordinates, $x=(a+b)/\sqrt2$, $y=(a-b)/\sqrt2$, and set $\xi=x^2/(1+s)$, $\eta=y^2/(1-s)$. Then $a^2-2sab+b^2=(1-s)x^2+(1+s)y^2$ and $(a-sb)(b-sa)=\tfrac12(1-s)^2x^2-\tfrac12(1+s)^2y^2$. Hence
\[
\partial_a\partial_b\varphi_2=\frac{\mathrm e^{-(\xi+\eta)/2}}{2\pi\sqrt q}\Bigl(\alpha_+ \xi-\alpha_- \eta+\gamma\Bigr),\qquad
\alpha_+=\frac{1}{2(1+s)},\ \alpha_-=\frac1{2(1-s)},\ \gamma=\frac{s}{q},
\]
and $(\xi,\eta)$ ranges over $[0,\infty)^2$. Write $f(\xi,\eta)=\mathrm e^{-(\xi+\eta)/2}(\alpha_+ \xi-\alpha_- \eta+\gamma)$; then $\sup|f|=\max\{\sup f,\ \sup(-f)\}$, and points where $f=0$ do not matter. \emph{Positive values.} If $f(\xi,\eta)>0$, then $\alpha_+ \xi+\gamma\ge\alpha_+ \xi-\alpha_- \eta+\gamma>0$ and $\mathrm e^{-\eta/2}\le1$, so $f(\xi,0)\ge f(\xi,\eta)$. Hence $\sup f=\sup_{\xi\ge0}f(\xi,0)=\sup_{\xi\ge0}\mathrm e^{-\xi/2}(\alpha_+ \xi+\gamma)$. The derivative in $\xi$ vanishes at $\xi^\ast=2-\gamma/\alpha_+=2(1-2s)/(1-s)$, which is admissible iff $s<\tfrac12$ and gives $2\alpha_+\,\mathrm e^{-\xi^\ast/2}=h_+(s)$. For $s\ge\frac12$ the function is decreasing on $[0,\infty)$, so the supremum is $f(0,0)=\gamma=h_+(s)$. \emph{Negative values.} If $f(\xi,\eta)<0$, then $\alpha_- \eta-\gamma\ge\alpha_- \eta-\alpha_+ \xi-\gamma>0$ and $\mathrm e^{-\xi/2}\le1$, so $-f(0,\eta)\ge-f(\xi,\eta)$. Hence $\sup(-f)=\sup_{\eta\ge0}\mathrm e^{-\eta/2}(\alpha_- \eta-\gamma)$, attained at $\eta^\ast=2+\gamma/\alpha_-=2(1+2s)/(1+s)$ with value $2\alpha_-\,\mathrm e^{-\eta^\ast/2}=\mathrm e^{-(1+2s)/(1+s)}/(1-s)$. Multiplying by $1/(2\pi\sqrt q)$ gives $h(s)$.

\emph{Monotonicity.} We have $2\pi h=\max\{g_+,g_-\}$ with
$g_-(s)=\mathrm e^{-(1+2s)/(1+s)}(1-s)^{-1}(1-s^2)^{-1/2}$ and $g_+(s)=h_+(s)(1-s^2)^{-1/2}$. A maximum of two nondecreasing functions is nondecreasing, so it suffices to treat $g_\pm$. Logarithmic differentiation gives, for $s\in[0,1)$,
\[
\frac{g_-'(s)}{g_-(s)}=\frac{1}{1-s}-\frac{1}{(1+s)^2}+\frac{s}{1-s^2}\ \ge\ 0,
\]
because $1/(1-s)\ge1\ge1/(1+s)^2$ and $s/(1-s^2)\ge0$. For $s\in[0,\frac12)$,
\[
\frac{g_+'(s)}{g_+(s)}=\frac{1}{(1-s)^2}-\frac{1}{1+s}+\frac{s}{1-s^2}\ \ge\ 0,
\]
because $1/(1-s)^2\ge1\ge1/(1+s)$ and $s/(1-s^2)\ge0$. For $s\in[\frac12,1)$, $g_+(s)=s(1-s^2)^{-3/2}$ is increasing. Finally both expressions for $h_+$ equal $\frac23$ at $s=\frac12$, so $g_+$ is continuous there and hence nondecreasing on all of $[0,1)$. Continuity of $h$ follows in the same way.
\end{Proof}

Define, for $\lambda\in[0,1)$,
\begin{equation}\label{eq:R}
R(\lambda)=\int_0^\lambda(\lambda-s)\,h(s)\,ds .
\end{equation}
Since $h(0)=\mathrm e^{-1}/(2\pi)$ and $h$ is continuous, $R(\lambda)=\lambda^2/(4\pi\mathrm e)+O(\lambda^3)$. Moreover $R(\lambda)/\lambda^2=\int_0^1(1-t)\,h(\lambda t)\,dt$, which is nondecreasing in $\lambda$ because $h$ is nondecreasing (Lemma~\ref{lem:h}). As $s\to1$, $h(s)\sim\frac{1}{2\pi}(1-s^2)^{-3/2}$, so $(\lambda-s)h(s)\le(1-s)h(s)$ is integrable and $R$ stays bounded: $R(1)=\lim_{\lambda\to1}R(\lambda)\approx0.1015$.
A cruder closed form follows from Mehler's expansion \citep[18.18.28]{DLMF} and Sz\'asz's inequality \citep{Szasz1951}, \citep[18.14.9]{DLMF}, which in the probabilists' normalisation reads $|He_n(x)|\varphi(x)\le\sqrt{n!/(2\pi)}$. Together they give $h(s)\le 1/(2\pi(1-s)^2)$ and therefore $R(\lambda)\le(-\lambda-\ln(1-\lambda))/(2\pi)\approx\lambda^2/(4\pi)$, which is $\mathrm e$ times larger for small $\lambda$. We know of no earlier closed form for the supremum in Lemma~\ref{lem:h}. Numerical values of $R$ are given in Table~\ref{tab:R}.

\begin{theorem}\label{thm:gauss}
For every $g:[r]\times[c]\to\R$ and $|\lambda|<1$,
\begin{equation}\label{eq:main}
\frac{\arcsin|\lambda|}{2\pi}\,\rho(g)\;\le\;S_\lambda(g)\;\le\;\frac{|\lambda|}{2\pi}\,\rho(g)+V(g)\,R(|\lambda|).
\end{equation}
The lower bound is attained by two-point marginals. Consequently
\[
S_\lambda(g)=\frac{|\lambda|}{2\pi}\,\rho(g)+O(\lambda^2),\qquad \lim_{\lambda\to0}\frac{S_\lambda(g)}{|\lambda|}=\frac{\rho(g)}{2\pi},
\]
where, by \eqref{eq:main}, the $O(\lambda^2)$ term $S_\lambda(g)-|\lambda|\rho(g)/(2\pi)$ itself (not merely an estimate of it) lies in the interval
\[
\Bigl[\tfrac{1}{2\pi}(\arcsin|\lambda|-|\lambda|)\,\rho(g),\ V(g)R(|\lambda|)\Bigr].
\] In particular, in the first-order sense as $\lambda\to0$, neither the constant $1/(2\pi)$ nor the index $\rho(g)$ can be replaced by anything smaller. Moreover,
\begin{equation}\label{eq:G}
S_\lambda(g)\le\rho(g)\,G(|\lambda|),\qquad G(\lambda)=\sum_{n\ge1}\frac{\lambda^nT_n^2}{4\,n!},\qquad T_n=\int_{\R}|\varphi^{(n)}(x)|\,dx\le\sqrt{n!},
\end{equation}
so that $G(\lambda)=\lambda/(2\pi)+G_2(\lambda)$ with $0\le G_2(\lambda)\le\lambda^2/(4(1-\lambda))$.
\end{theorem}

\begin{Proof}
\emph{Representation.} Let $a_i=\Phi^{-1}(u_i)$ and $b_j=\Phi^{-1}(v_j)$, with values in $[-\infty,\infty]$. The thresholds are nondecreasing in $i$ and in $j$. A cell with an infinite threshold contributes nothing to \eqref{eq:identity} (Remark~\ref{rem:inf}), so in the next display we may assume $a_i,b_j$ finite. By Lemma~\ref{lem:plackett},
\[
\begin{aligned}
c_{ij}:=C^{\Ga}_\lambda(u_i,v_j)-u_iv_j&=\int_0^\lambda\varphi_2(a_i,b_j;s)\,ds\\
&=\lambda\varphi(a_i)\varphi(b_j)+\int_0^\lambda(\lambda-s)\,\partial_a\partial_b\varphi_2(a_i,b_j;s)\,ds ,
\end{aligned}
\]
where the last step is Taylor's formula with integral remainder for $\lambda\mapsto\Phi_2(a_i,b_j;\lambda)$ around $\lambda=0$, using $\varphi_2(a,b;0)=\varphi(a)\varphi(b)$ and the second identity of Lemma~\ref{lem:plackett}. By Lemma~\ref{lem:h} the remainder is bounded by $R(|\lambda|)$ in absolute value, uniformly in $(a_i,b_j)$.

\emph{Main term.} Put $w_i=\varphi(a_i)/\varphi(0)$ and $z_j=\varphi(b_j)/\varphi(0)$. Since $\varphi$ is unimodal and the thresholds are nondecreasing, $w$ and $z$ are unimodal sequences in $[0,1]$; infinite thresholds give $w_i=0$ or $z_j=0$ and do not break unimodality (Remark~\ref{rem:inf}). By Lemma~\ref{lem:layer},
\[
\lvert \lambda\sum_{i,j}\Delta_{ij}(g)\varphi(a_i)\varphi(b_j)\rvert =\frac{|\lambda|}{2\pi}\,|w^\top\Delta(g)z|\le\frac{|\lambda|}{2\pi}\,\rho(g).
\]
Together with $|\sum_{ij}\Delta_{ij}(g)(c_{ij}-\lambda\varphi(a_i)\varphi(b_j))|\le V(g)R(|\lambda|)$ this gives the upper bound in \eqref{eq:main} for every pair of marginals.

\emph{Lower bound.} Let $I=\{i_1,\dots,i_2\}$ and $J=\{j_1,\dots,j_2\}$ attain the maximum in \eqref{eq:rho}. Let $X$ take the values $i_1$ and $i_2+1$ with probability $\frac12$ each, and let $Y$ take $j_1$ and $j_2+1$ likewise. Then $u_i=\frac12$ for $i\in I$, while $u_i\in\{0,1\}$ for $i\notin I$, and similarly for $v$. By \eqref{eq:identity} and Sheppard's formula $\Phi_2(0,0;\lambda)=\frac14+\arcsin(\lambda)/(2\pi)$ \citep{Sheppard1899},
\[
B=\Bigl(\sum_{i\in I,j\in J}\Delta_{ij}(g)\Bigr)\bigl(\Phi_2(0,0;\lambda)-\tfrac14\bigr)=\pm\rho(g)\,\frac{\arcsin\lambda}{2\pi}.
\]
The asymptotic statements follow since $\arcsin|\lambda|=|\lambda|+O(\lambda^3)$ and $R(|\lambda|)=O(\lambda^2)$.

\emph{Tetrachoric bound.} By Lemma~\ref{lem:plackett}, $\partial_s^k\varphi_2(a,b;s)|_{s=0}=\varphi^{(k)}(a)\varphi^{(k)}(b)$, which gives the tetrachoric series \citep{Pearson1900}
\[
c_{ij}=\sum_{n\ge1}\frac{\lambda^n}{n!}\,\varphi^{(n-1)}(a_i)\,\varphi^{(n-1)}(b_j),\qquad |\lambda|<1,
\]
equivalently Mehler's expansion integrated over the quadrant. Since $\varphi^{(n)}=(-1)^n\mathrm{He}_n\varphi$, we have $T_n=\E|\mathrm{He}_n(Z)|\le(\E\,\mathrm{He}_n(Z)^2)^{1/2}=\sqrt{n!}$ for standard normal $Z$, by the Cauchy--Schwarz inequality. The thresholds $a_i$ are nondecreasing and $\varphi^{(n-1)}$ vanishes at $\pm\infty$, so the zero-padded sequence $(\varphi^{(n-1)}(a_i))_i$ has total variation at most that of $\varphi^{(n-1)}$ on $\R$, which is $T_n$; thus $\nu\bigl((\varphi^{(n-1)}(a_i))_i\bigr)\le T_n/2$, and $|\varphi^{(n-1)}|\le T_n/2$ pointwise. In particular the series converges absolutely, its $n$th term being at most $|\lambda|^nT_n^2/(4\,n!)\le|\lambda|^n/4$. By Lemma~\ref{lem:layer}, the $n$th term of $\sum_{ij}\Delta_{ij}(g)c_{ij}$ is at most $\rho(g)|\lambda|^nT_n^2/(4\,n!)$ in absolute value. For $n=1$, $T_1=2\varphi(0)$ and this is $|\lambda|\rho(g)/(2\pi)$. Summing over $n$ gives \eqref{eq:G}, and $T_n^2\le n!$ gives the bound on $G_2$.
\end{Proof}

Both remainder bounds hold, so
\[
S_\lambda(g)\le\frac{|\lambda|}{2\pi}\rho(g)+\min\bigl\{V(g)R(|\lambda|),\ \rho(g)G_2(|\lambda|)\bigr\}.
\]
The constants $T_n$ are explicit: the extrema of $\varphi^{(n-1)}$ lie at the zeros $x_1<\dots<x_n$ of $\mathrm{He}_n$ and alternate in sign, so $T_n=2\sum_k|\mathrm{He}_{n-1}(x_k)|\varphi(x_k)$; for instance $T_1=2\varphi(0)$ and $T_2=4\varphi(1)$. Hence $G_2(\lambda)=0.117\lambda^2+O(\lambda^3)$, against $R(\lambda)=0.029\lambda^2+O(\lambda^3)$, so the tetrachoric bound is the sharper one when $V(g)/\rho(g)$ exceeds about $4$ and $\lambda$ is small. Numerically $G_2(0.2)\approx0.0056$, $G_2(0.3)\approx0.0140$, $G_2(0.6)\approx0.0860$ and $G_2(0.9)\approx0.509$. The bound $\rho(g)G(|\lambda|)$ depends on $g$ only through $\rho(g)$: at every $|\lambda|<1$, $\rho(g)$ determines $S_\lambda(g)$ up to the factor $2\pi G(|\lambda|)/\arcsin|\lambda|$, which equals $1.08$, $1.17$, $1.27$, $1.77$ and $3.66$ at $|\lambda|=0.1$, $0.2$, $0.3$, $0.6$ and $0.9$. For the individual terms of the series the bound of Lemma~\ref{lem:layer} is exact for $n=1$; in numerical maximisation over random tables and thresholds the largest values found were $83\%$, $55\%$, $46\%$, $39\%$ and $41\%$ of $T_n^2/4$ for $n=2,\dots,6$.

\begin{remark}[Infinite thresholds and attainment]\label{rem:inf}
Thresholds equal to $\pm\infty$ occur exactly when some $u_i$ or $v_j$ lies in $\{0,1\}$. With the conventions $\varphi(\pm\infty)=0$, $\Phi_2(-\infty,b)=0$ and $\Phi_2(+\infty,b)=\Phi(b)$, such cells have $c_{ij}=0$. The sequence $w$ stays unimodal, since repeated thresholds (grades of probability zero) are allowed. The two-point marginals in the proof produce thresholds $-\infty,0,+\infty$ exactly, so the lower bound is attained and not merely approached.
\end{remark}

When $\Delta(g)$ has constant sign the problem is solved exactly for every $\lambda$, not only to first order. This follows from the elementary bound $0\le c_{ij}\le\arcsin(\lambda)/(2\pi)$ for $\lambda\ge0$, which holds because $0<\varphi_2(a,b;s)\le\varphi_2(0,0;s)=1/(2\pi\sqrt{1-s^2})$. This is an upper bound of the Vitali type considered by \citet{Beare2009}, split into positive and negative parts.

\begin{proposition}\label{prop:sign}
Let $P(g)$ and $Q(g)$ be the sums of the positive parts and of the negative parts of the entries of $\Delta(g)$, so that $V=P+Q$. Then $S_\lambda(g)\le\frac{\arcsin|\lambda|}{2\pi}\max\{P(g),Q(g)\}$. If $\Delta(g)\ge0$ or $\Delta(g)\le0$ entrywise, then $\rho(g)=V(g)=|N(g)|$ and
\[
S_\lambda(g)=\frac{\arcsin|\lambda|}{2\pi}\,|N(g)|\qquad\text{for all }|\lambda|<1.
\]
\end{proposition}

\begin{Proof}
For $\lambda\ge0$ every $c_{ij}$ lies in $[0,\arcsin(\lambda)/(2\pi)]$, so $-Q\arcsin(\lambda)/(2\pi)\le B\le P\arcsin(\lambda)/(2\pi)$. The case $\lambda<0$ is symmetric. If $\Delta\ge0$, the sum over a rectangle can only grow when the rectangle is enlarged, so the maximum in \eqref{eq:rho} is attained by the full grid and $\rho=\sum_{i,j}\Delta_{ij}=N=P=V$ (similarly for $\Delta\le0$). The upper bound then coincides with the lower bound of Theorem~\ref{thm:gauss}.
\end{Proof}

Combining Theorem~\ref{thm:gauss} with Proposition~\ref{prop:sign} gives
\[
S_\lambda(g)\le\min\Bigl\{\frac{|\lambda|}{2\pi}\rho+\min\{VR(|\lambda|),\rho G_2(|\lambda|)\},\ \frac{\arcsin|\lambda|}{2\pi}\max(P,Q)\Bigr\};
\]
the second bound is the better one when $\max(P,Q)$ is close to $\rho(g)$ or $\lambda$ is close to $1$. When mixed differences of both signs cancel, Theorem~\ref{thm:gauss} is sharper; for the matrix $M$ of Example~\ref{ex:matrix}, $\max(P,Q)=3\rho(M)$.

For a supermodular table ($\Delta\ge0$, ``and''-like aggregation) and positively dependent inputs, the independence computation therefore underestimates $\E g$ by up to $|N(g)|\arcsin(\lambda)/(2\pi)$, and this bound is sharp. For a mixed-sign table Theorem~\ref{thm:gauss} replaces $|N|$ by $\rho\ge|N|$ to first order.

\section{OTHER COPULA FAMILIES}\label{sec:families}

The proof of Theorem~\ref{thm:gauss} used only two features of the Gaussian copula: the first-order departure from independence is a product $A(u)A(v)$, and $A$ is unimodal. We now isolate these features, allowing non-exchangeable products $A(u)B(v)$.

\begin{theorem}\label{thm:families}
Let $\{C_\theta\}$ be a family of bivariate copulas with $C_{\theta_0}=\Pi$. Let $\Theta$ be an interval of parameters containing $\theta_0$ and with $\theta_0$ not isolated in $\Theta$; one-sided neighbourhoods such as $\Theta=[\theta_0,\theta_0+\delta]$ are allowed. Suppose there exist $L\ge0$, a sign $\sigma\in\{-1,1\}$ and unimodal functions $A,B:[0,1]\to[0,\infty)$ (that is, $A$ is nondecreasing on $[0,u^\ast]$ and nonincreasing on $[u^\ast,1]$ for some $u^\ast$, and likewise $B$ with some $v^\ast$) such that
\begin{equation}\label{eq:rankone}
\lvert C_\theta(u,v)-uv-\sigma(\theta-\theta_0)A(u)B(v)\rvert \le L(\theta-\theta_0)^2
\end{equation}
for all $u,v\in[0,1]$ and $\theta\in\Theta$.
Let $\kappa=\max_uA(u)\cdot\max_vB(v)=A(u^\ast)B(v^\ast)$. Then for every $g$ and $\theta\in\Theta$,
\[
|\theta-\theta_0|\,\kappa\,\rho(g)-L\rho(g)(\theta-\theta_0)^2\;\le\;S_{C_\theta}(g)\;\le\;|\theta-\theta_0|\,\kappa\,\rho(g)+L\,V(g)(\theta-\theta_0)^2 .
\]
\end{theorem}

\begin{Proof}
Since $C_\theta(0,v)=0$, putting $u=0$ in \eqref{eq:rankone} gives $|\theta-\theta_0|\,A(0)B(v)\le L(\theta-\theta_0)^2$; dividing by $|\theta-\theta_0|$ and letting $\theta\to\theta_0$ within $\Theta$ yields $A(0)B(v)=0$ for all $v$. The same argument with $C_\theta(1,v)=v$, and with the roles of $u$ and $v$ exchanged, yields $A(1)B(v)=0$ and $A(u)B(0)=A(u)B(1)=0$. We may assume $A,B\not\equiv0$, so $A$ and $B$ vanish at $0$ and $1$. Consequently $u^\ast,v^\ast\in(0,1)$ and $\kappa>0$. (If $A\equiv0$ or $B\equiv0$, the first-order term vanishes and the theorem only gives $S_{C_\theta}(g)\le LV(g)(\theta-\theta_0)^2$.) \emph{Upper bound.} Write $\mathcal B$ for the bias. By \eqref{eq:identity}, $\mathcal B=\sigma(\theta-\theta_0)\sum_{ij}\Delta_{ij}A(u_i)B(v_j)+\sum_{ij}\Delta_{ij}e_{ij}$ with $|e_{ij}|\le L(\theta-\theta_0)^2$. The sequences $A(u_i)/A(u^\ast)$ and $B(v_j)/B(v^\ast)$ are unimodal and take values in $[0,1]$, because $u_i$ and $v_j$ are nondecreasing and $A,B$ are unimodal. Lemma~\ref{lem:layer} bounds the first sum by $\kappa\rho(g)$. \emph{Lower bound.} With $I,J$ as in the proof of Theorem~\ref{thm:gauss}, let $X=i_1$ with probability $u^\ast$ and $X=i_2+1$ with probability $1-u^\ast$, and let $Y=j_1$ with probability $v^\ast$ and $Y=j_2+1$ with probability $1-v^\ast$. Then $u_i=u^\ast$ on $I$ and $u_i\in\{0,1\}$ elsewhere, and $v_j=v^\ast$ on $J$ and $v_j\in\{0,1\}$ elsewhere. Hence
$\mathcal B=\bigl(\sum_{I\times J}\Delta_{ij}\bigr)\bigl[\sigma(\theta-\theta_0)\kappa+e\bigr]$ with $|e|\le L(\theta-\theta_0)^2$. Taking absolute values, $|\mathcal B|\ge\rho(g)\bigl(|\theta-\theta_0|\kappa-L(\theta-\theta_0)^2\bigr)$ whatever the sign of $\sum_{I\times J}\Delta_{ij}$.
\end{Proof}

Unimodality is needed only to make the two bounds coincide: without it the lower bound still holds by Proposition~\ref{prop:general}, since $\varepsilon(C_\theta)\ge|\theta-\theta_0|\sup A\sup B-L(\theta-\theta_0)^2$, while the upper bound changes as follows. If $A,B$ are merely of bounded variation (they still vanish at $0$ and $1$ by the first step of the proof), then for nondecreasing $u_i$ the zero-padded sequence $(A(u_i))_i$ satisfies $\nu\le\mathrm{TV}(A)/2$, where $\mathrm{TV}$ is the total variation on $[0,1]$, and Lemma~\ref{lem:layer} gives the upper bound with $\kappa$ replaced by $\mathrm{TV}(A)\,\mathrm{TV}(B)/4$. For unimodal $A,B$ this equals $A(u^\ast)B(v^\ast)=\kappa$.

For the Gaussian family $A(u)=\varphi(\Phi^{-1}(u))$, $u^\ast=\frac12$, $\kappa=1/(2\pi)$, and \eqref{eq:rankone} holds with $L=\sup_{|\lambda|\le\delta}R(|\lambda|)/\lambda^2=R(\delta)/\delta^2$ by the symmetry $h(-s)=h(s)$ and the monotonicity of $R(\lambda)/\lambda^2$ noted after \eqref{eq:R}. For $\delta=0.6$ this gives $L\in[0.034956,\,0.034958]$: since $h$ is nondecreasing, $R(\delta)$ lies between the lower and upper Riemann sums $\sum_kh(s_k)\int_{s_k}^{s_{k+1}}(\delta-s)\,ds$ and $\sum_kh(s_{k+1})\int_{s_k}^{s_{k+1}}(\delta-s)\,ds$, which we evaluated in Arb on $2\cdot10^4$ equal subintervals. Theorem~\ref{thm:gauss} is the sharper special case. Table~\ref{tab:families} lists other common one-parameter families. Their first derivatives at independence were obtained symbolically and are exact. For FGM the expansion is exact ($L=0$), so $S_{C_\theta}(g)=|\theta|\rho(g)/16$ for all $\theta\in[-1,1]$. For AMH, $C_\theta-uv-\theta A(u)A(v)=\theta^2uv\,w^2/(1-\theta w)$ with $w=(1-u)(1-v)$, and $\max_{u,v}uv\,w^2=(4/27)^2$ gives $L=16/(729(1-\delta))$. This bound blows up as $\delta\to1$, but $L$ itself stays finite on the whole range $|\theta|\le1$: the remainder divided by $\theta^2$, $uv\,w^2/(1-\theta w)$, is increasing in $\theta$ for $w\ge0$, so $L=\sup_{u,v}uv\,w^2/(u+v-uv)$, the value at $\theta=1$. Interval branch and bound in Arb (with $uv/(u+v-uv)\le\min(u,v)$ near the axes) gives $L\in[0.045749,\,0.045754]$, attained near $u=v=0.22$. For these three families (Gaussian, FGM, AMH) Theorem~\ref{thm:families} therefore applies with the stated constants. For Frank, Clayton and Plackett the first derivative, and hence $A$ and $\kappa$, are exact. The existence of a finite $L$ can be shown for Frank and Plackett as follows. For Frank, the argument of the logarithm in $C_\theta$ equals $1-\theta uv\,k$ with $k=e_1(\theta u)e_1(\theta v)/e_1(\theta)$ and $e_1(y)=(1-\mathrm e^{-y})/y$, which is analytic and positive on $\R$; the argument also equals $\mathrm e^{-\theta C_\theta}\ge\mathrm e^{-|\theta|}$. For Plackett, $C_\theta=2\theta uv/(S+\sqrt D)$ with $S=1+(\theta-1)(u+v)\ge0$ and $D=1+2(\theta-1)(u+v-2uv)+(\theta-1)^2(u-v)^2$, which is positive on $[0,1]^2\times[\frac12,\frac32]$. In both cases $C_\theta(u,v)$ extends to a real-analytic function of $(u,v,\theta)$ on a neighbourhood of the compact set $[0,1]^2\times[\theta_0-\delta,\theta_0+\delta]$, so $\partial_\theta^2C_\theta$ is bounded there and \eqref{eq:rankone} holds with $L=\frac12\sup|\partial_\theta^2C_\theta|<\infty$. Theorem~\ref{thm:families} therefore applies. Explicit values of $L$ were obtained by a computer-assisted proof, described at the end of this section. For Clayton we set $C_0=\Pi$, the limit of $C_\theta$ as $\theta\downarrow0$, and take $\Theta=[0,\infty)$ (positive dependence; for $\theta<0$ the copula has a zero region and is not smooth in $\theta$). We obtain an explicit, if crude, certified constant. Put $x=-\ln u$, $y=-\ln v$, $m=\max(x,y)$ and $W(\theta)=\mathrm e^{\theta x}+\mathrm e^{\theta y}-1\ge\max(\mathrm e^{\theta x},\mathrm e^{\theta y})$. Then $\ln C_\theta=-\Psi(\theta)/\theta$ with $\Psi=\ln W$ and $\Psi(0)=0$, so $\partial_\theta\ln C_\theta=-\int_0^1t\,\Psi''(t\theta)\,dt$ and $\partial^2_\theta\ln C_\theta=-\int_0^1t^2\Psi'''(t\theta)\,dt$. The bound on $W$ gives $0\le W^{(k)}/W\le x^k+y^k$, hence $|\Psi''|\le(x+y)^2\le4m^2$ and $|\Psi'''|\le30m^3$. With $C_\theta\le\min(u,v)=\mathrm e^{-m}$ this yields
$|\partial_\theta^2C_\theta|=C_\theta\,|(\partial_\theta\ln C_\theta)^2+\partial_\theta^2\ln C_\theta|\le\mathrm e^{-m}(4m^4+10m^3)\le31.2$
for all $\theta>0$ and $u,v\in(0,1]$. Since $C_\theta$ and $\partial_\theta C_\theta$ extend continuously to $\theta=0$, with $\partial_\theta C_\theta\to uv\ln u\ln v$, Taylor's formula on $[0,\theta]$ shows that \eqref{eq:rankone} holds with $L\le15.6$ on $\Theta=[0,\infty)$, without computer assistance. Much better: on $\Theta=[0,1]$ the computer-assisted bound below gives $L\le0.0691$, and for $\theta\ge1$ trivially $|C_\theta-uv-\theta A(u)A(v)|\le\frac14+\theta\mathrm e^{-2}\le(\frac14+\mathrm e^{-2})\theta^2$, since $|C-\Pi|\le\frac14$ for every copula. Hence $L\le\frac14+\mathrm e^{-2}<0.386$ on $\Theta=[0,\infty)$. The same argument gives, for any family, $L\le1/(4\delta^2)+\kappa/\delta$ on $\{|\theta-\theta_0|\ge\delta\}$. In numerical checks at $|\theta-\theta_0|=0.02$ the ratio $S_{C_\theta}(g)/(|\theta-\theta_0|\kappa\rho(g))$ lay between $0.986$ and $1.005$ for all five non-Gaussian families. In all families of Table~\ref{tab:families} $B=A$, $A\ge0$, $A(0)=A(1)=0$ and $A$ is unimodal with an interior maximiser $u^\ast$, which can be checked directly from the formulas.

\textbf{FGM-type families.} If $C_\theta(u,v)=uv+\theta A(u)B(v)$ exactly, then $L=0$ and $S_{C_\theta}(g)=|\theta|\kappa\rho(g)$ for every admissible $\theta$. This covers the generalised FGM (Sarmanov-type) families with nonnegative unimodal $A,B$, including the semiparametric families $uv+\theta A(u)A(v)$ of \citet{AmblardGirard2002} and $uv+\theta f(u)g(v)$ of \citet{RodriguezLallenaUbedaFlores2004}, and the single-term cases of the perturbations of independence $uv+\sum_k\lambda_k\Phi_k(u)\Phi_k(v)$ studied by \citet{Longla2024}; for several terms Lemma~\ref{lem:layer} still gives an upper bound, as in the Remark below. For example, the Huang--Kotz copula $uv[1+\theta(1-u^p)(1-v^p)]$, $p>0$, has $A(u)=u(1-u^p)$, $u^\ast=(1+p)^{-1/p}$ and $\kappa=p^2(1+p)^{-2-2/p}$ \citep{HuangKotz1999}. Its density $1+\theta(1-(1+p)u^p)(1-(1+p)v^p)$ is nonnegative exactly for $-\min(1,p^{-2})\le\theta\le1/p$.

\textbf{Archimedean families.} Let $C_\theta(u,v)=\psi_\theta^{-1}(\psi_\theta(u)+\psi_\theta(v))$ with generators normalised so that $\psi_{\theta_0}(t)=-\ln t$. Assume $(\theta,t)\mapsto\psi_\theta(t)$ is continuously differentiable on $\Theta\times(0,1]$ with $\partial_t\psi_\theta<0$. Put
\[
\chi(x)=\partial_\theta\psi_\theta(\mathrm e^{-x})\big|_{\theta=\theta_0},\qquad x\ge0 .
\]
Then $\chi(0)=0$, because $\psi_\theta(1)=0$ for all $\theta$. Implicit differentiation of $\psi_\theta(C_\theta)=\psi_\theta(u)+\psi_\theta(v)$ at $\theta_0$ (the bivariate case of \citealp[Prop.~3.2]{GenestQuessyRemillard2007}, who computed these derivatives for Frank, AMH, Clayton, Gumbel--Barnett and Gumbel in a study of local power of tests of independence), using $\partial_t\psi_{\theta_0}(t)=-1/t$ and $C_{\theta_0}=uv$, gives
\begin{equation}\label{eq:archfirst}
\partial_\theta C_\theta(u,v)\big|_{\theta_0}=uv\,\bigl[\chi(x+y)-\chi(x)-\chi(y)\bigr],\qquad x=-\ln u,\ y=-\ln v,
\end{equation}
that is, $uv$ times the Cauchy difference of $\chi$.

\begin{proposition}\label{prop:arch}
Let $\chi$ be continuous. The first-order term \eqref{eq:archfirst} has the form $\sigma A(u)A(v)$ with $\sigma=\pm1$, $A\ge0$ and $A\not\equiv0$ if and only if, up to adding a linear function $bx$ to $\chi$ (which corresponds to rescaling the generators), for some $c>0$ either
(i) $\chi(x)=\sigma c^2x^2/2$, and then $A(u)=c\,u\ln(1/u)$; or
(ii) $\chi(x)=\sigma c^2\beta^{-2}(\mathrm e^{\beta x}-1-\beta x)$ with $\beta\ne0$, and then $A(u)=|c/\beta|\,u\,|u^{-\beta}-1|$.
In case (i) $A$ is concave. If, in addition, the remainder condition \eqref{eq:rankone} holds, then in case (ii) necessarily $\beta<1$ and $A$ is concave. In both cases $A$ is then unimodal, so Theorem~\ref{thm:families} applies.
\end{proposition}

\begin{Proof}
Write $\chi(x+y)-\chi(x)-\chi(y)=\sigma F(x)F(y)$ with $F(x)=A(\mathrm e^{-x})\mathrm e^{x}\ge0$, $x\ge0$. Fixing $y$ with $F(y)\ne0$ shows that $F$ is continuous, since $\chi$ is. Putting $y=0$ and using $\chi(0)=0$ gives $F(0)F(x)=0$, so $F(0)=0$. The Cauchy difference $G(x,y)=F(x)F(y)$ satisfies the cocycle identity $G(x+y,z)+G(x,y)=G(x,y+z)+G(y,z)$, which can be checked by expanding $\chi$. Hence $F(z)[F(x+y)-F(y)]=F(x)[F(y+z)-F(y)]$. Fix $z$ with $F(z)\ne0$; then $F(x+y)=F(y)+F(x)H(y)$ with $H(y)=[F(y+z)-F(y)]/F(z)$. Exchanging $x$ and $y$ gives $F(x)(H(y)-1)=F(y)(H(x)-1)$, so $H-1=kF$ for a constant $k$, and
$F(x+y)=F(x)+F(y)+kF(x)F(y)$.
If $k=0$, $F$ is a continuous additive function on $[0,\infty)$, hence linear, so $F(x)=cx$, with $c>0$ because $F\ge0$ and $F\not\equiv0$; the Cauchy difference of $\sigma c^2x^2/2$ is $\sigma c^2xy=\sigma F(x)F(y)$. If $k\ne0$, $E=1+kF$ is continuous with $E(x+y)=E(x)E(y)$ for $x,y\ge0$ and $E(0)=1$. The equation holds only on the half-line, but this suffices: $E(x)=E(x/2)^2\ge0$, and a zero $E(x_0)=0$ would give $E(x_0/2^n)=0$ for all $n$ and, by continuity, $E(0)=0$; so $E>0$, $\ln E$ is continuous and additive on $[0,\infty)$, and $E(x)=\mathrm e^{\beta x}$ with $\beta\ne0$ \citep{Aczel1966} and $F(x)=(\mathrm e^{\beta x}-1)/k$; the Cauchy difference of $\sigma k^{-2}(\mathrm e^{\beta x}-1-\beta x)$ is $\sigma k^{-2}(\mathrm e^{\beta x}-1)(\mathrm e^{\beta y}-1)=\sigma F(x)F(y)$, which is case (ii) with $c=|\beta/k|$. In each case the stated $\chi$ is a particular solution, and the general solution differs from it by a solution of $\chi(x+y)=\chi(x)+\chi(y)$, that is, by $bx$, because continuous additive functions are linear. The formulas for $A(u)=u\,|F(-\ln u)|$ follow. The converse is a direct computation. Now assume \eqref{eq:rankone}. As in the proof of Theorem~\ref{thm:families}, $A(0)=0$. Moreover $A$ is continuous: copulas are $1$-Lipschitz in each argument, so \eqref{eq:rankone} gives $t|A(u)-A(u')|A(v)\le2|u-u'|+2Lt^2$ for $t=|\theta-\theta_0|$, and letting $u'\to u$ and then $t\to0$ (along $\Theta$) gives continuity; here $A(v)>0$ for some $v$. (The boundary identity $C_\theta(0,v)=0$ alone would not do, since it says nothing about $\partial_\theta C_\theta(u,v)$ as $u\to0$.) Hence $A(u)\to0$ as $u\to0$, and case (ii) requires $u^{1-\beta}\to0$ as $u\to0$ (for $\beta>0$, $A(u)\sim|c/\beta|u^{1-\beta}$; for $\beta<0$, $A(u)\sim|c/\beta|u$), that is, $\beta<1$. Concavity follows from $A''(u)=-c/u$ in case (i), $A''(u)\propto-\beta(1-\beta)u^{-\beta-1}$ for $0<\beta<1$, and $A''(u)\propto-p(1+p)u^{p-1}$ for $\beta=-p<0$.
\end{Proof}

For AMH and Frank, $\chi(x)=\mathrm e^{-x}-1$ and $\chi(x)=(\mathrm e^{-x}-1)/2$ are of type (ii) with $\beta=-1$. For Clayton, $\chi(x)=x^2/2$ is of type (i). The Gumbel--Barnett family $C_\theta(u,v)=uv\exp(-\theta\ln u\ln v)$, $0\le\theta\le1$ \citep[family 4.2.9]{Nelsen2006}, has $\chi(x)=-x^2/2$, type (i) with $\sigma=-1$. For it the remainder is explicit: $C_\theta-uv+\theta uv\ln u\ln v=uv(\mathrm e^{-\theta xy}-1+\theta xy)\in[0,\theta^2(u\ln^2u)(v\ln^2v)/2]$, so $L=\frac12(4\mathrm e^{-2})^2=8\mathrm e^{-4}\approx0.147$. By contrast, the Gumbel family ($\theta_0=1$) has $\chi(x)=x\ln x$, and the Joe family, $\psi_\theta(t)=-\ln(1-(1-t)^\theta)$ with $\theta_0=1$, has $\chi(x)=(\mathrm e^{x}-1)\ln(1-\mathrm e^{-x})$. Neither is of type (i) or (ii), so by Proposition~\ref{prop:arch} their first-order terms are not of product form. We have checked \eqref{eq:archfirst} numerically against finite differences for all six Archimedean families mentioned. The product form of \eqref{eq:archfirst} is equivalent to the factorisation $\partial_\theta c_\theta|_{\theta_0}=\sigma A'(u)A'(v)$ of the score of the copula density, that is, to the locally most powerful rank test of independence against the family being a linear rank statistic \citep{GenestVerret2005}: Spearman's statistic for AMH and Frank and Savage-type scores for Clayton. We have not found the characterisation of Proposition~\ref{prop:arch} in that literature.

\textbf{Certified constants.} For Clayton, Plackett and Frank we bounded $\frac12\sup|\partial_\theta^2C_\theta|$ rigorously by interval (ball) arithmetic in Arb \citep{Johansson2017}, through python-flint, with branch and bound over boxes in $(u,v,\theta)$; the boxes cover the domain exactly, every series is truncated with an explicit tail bound, and a box is bisected until the upper end of its enclosure is within a factor $1.02$ of a certified lower bound for the supremum. For Clayton the formulas are rewritten through the cumulants of $\Psi=\ln W$ to avoid the removable singularity at $\theta=0$, the tail $\max(x,y)>30$ being covered by the analytic bound above; for Frank a centred form with forward-mode differentiation on balls is used, and the symmetries of the family reduce the domain to $\theta\in[0,1]$, $u\le v$, $u+v\le1$. The certified enclosures of $\frac12\sup|\partial_\theta^2C_\theta|$ are $[0.06766,\,0.06903]$ for Clayton on $[0,1]$ (the lower end is $\frac12\mathrm e^{-2}=0.0676676\ldots$ rounded down, the limit as $\theta\downarrow0$ at $u=v=1/\mathrm e$), $[0.11091,\,0.11314]$ for Plackett on $[\frac12,\frac32]$ and $[0.001368,\,0.001396]$ for Frank on $[-1,1]$. These are admissible values of $L$ in \eqref{eq:rankone} but not necessarily the smallest ones (on a grid the smallest are about $0.067$, $0.047$ and $0.00096$); a smaller $L$ only widens the range $|\theta-\theta_0|<\kappa/L$ in which the lower bound of Theorem~\ref{thm:families} is informative. The code is available with the paper.

\begin{table}[t]
\centering
\footnotesize
\setlength{\tabcolsep}{2pt}
\caption{First-order structure of common copula families near independence. Parameterisations follow \citet{Nelsen2006}; $\kappa$ depends on the parameterisation, $\rho(g)$ does not. All constants are proved: those marked $^\ast$ by a computer-assisted (interval-arithmetic) argument, see text, the others analytically. The $^\ast$ values bound $\frac12\sup|\partial_\theta^2C_\theta|$ and may exceed the smallest admissible $L$ in \eqref{eq:rankone}.}\label{tab:families}
\resizebox{\textwidth}{!}{%
\begin{tabular}{@{}lllccl@{}}
\toprule
Family & $\partial_\theta C_\theta|_{\theta_0}$ & $A(u)$ & $u^\ast$ & $\kappa$ & $L$ (range of $\theta$)\\
\midrule
Gaussian ($\theta_0=0$) & $\varphi(\Phi^{-1}u)\varphi(\Phi^{-1}v)$ & $\varphi(\Phi^{-1}u)$ & $1/2$ & $1/(2\pi)$ & $0.0350^\ast$ ($|\lambda|\le0.6$)\\
FGM ($\theta_0=0$) & $uv(1-u)(1-v)$ & $u(1-u)$ & $1/2$ & $1/16$ & $0$ ($|\theta|\le1$)\\
AMH ($\theta_0=0$) & $uv(1-u)(1-v)$ & $u(1-u)$ & $1/2$ & $1/16$ & $16/[729(1-\delta)]$ ($|\theta|\le\delta$); $0.0458^\ast$ ($|\theta|\le1$)\\
Frank ($\theta_0=0$) & $\tfrac12uv(1-u)(1-v)$ & $u(1-u)/\sqrt2$ & $1/2$ & $1/32$ & $0.00140^\ast$ ($|\theta|\le1$)\\
Clayton ($\theta_0=0$) & $uv\ln u\ln v$ & $u\ln(1/u)$ & $1/\mathrm e$ & $\mathrm e^{-2}$ & $0.0691^\ast$ ($0\le\theta\le1$); $0.386$ ($\theta\ge0$)\\
Plackett ($\theta_0=1$) & $uv(1-u)(1-v)$ & $u(1-u)$ & $1/2$ & $1/16$ & $0.1132^\ast$ ($|\theta-1|\le\frac12$)\\
Huang--Kotz ($\theta_0=0$) & $uv(1-u^p)(1-v^p)$ & $u(1-u^p)$ & $(1+p)^{-1/p}$ & $p^2(1+p)^{-2-2/p}$ & $0$ (admissible $\theta$)\\
Gumbel--Barnett ($\theta_0=0$) & $-uv\ln u\ln v$ & $u\ln(1/u)$ & $1/\mathrm e$ & $\mathrm e^{-2}$ & $8\mathrm e^{-4}$ ($0\le\theta\le1$)\\
\bottomrule
\end{tabular}}
\end{table}

\textbf{Range of the lower bound.} The lower bound of Theorem~\ref{thm:families}, $|\theta-\theta_0|\kappa\rho(g)-L\rho(g)(\theta-\theta_0)^2$, is positive only for $|\theta-\theta_0|<\kappa/L$. With the constants of Table~\ref{tab:families}, $\kappa/L$ is $4.55$ for the Gaussian family ($|\lambda|\le0.6$), $1.37$ for AMH ($|\theta|\le1$), $22.4$ for Frank, $1.96$ for Clayton on $[0,1]$ and $0.35$ with $L=0.386$ on $[0,\infty)$, $0.92$ for Gumbel--Barnett and $0.55$ for Plackett. The lower bound is therefore informative on the whole stated range of $\theta$, except for Clayton with $\theta>1$, where only $L=0.386$ is available, and Gumbel--Barnett with $\theta>0.92$; for FGM-type families ($L=0$) it is exact.

\begin{remark}[A family outside the scope]
For the Gumbel copula $C_\theta(u,v)=\exp\{-(x^\theta+y^\theta)^{1/\theta}\}$ with $x=-\ln u$, $y=-\ln v$, one finds
$\partial_\theta C_\theta|_{\theta=1}=uv\,[(x+y)\ln(x+y)-x\ln x-y\ln y]$.
This is not of product form (Proposition~\ref{prop:arch}), so Lemma~\ref{lem:layer} does not apply directly. In fact $\rho(g)$ no longer determines the first-order worst case exactly, although it does up to a factor $1.25$. Write $T=\partial_\theta C_\theta|_{\theta=1}$; the first-order coefficient is $\sup|\sum_{ij}\Delta_{ij}T(u_i,v_j)|$, the supremum being over all nondecreasing $(u_i)$, $(v_j)$. For $g=\1\{M=4\}$ it equals $\max_{u,v}T(u,v)$. Interval branch and bound in Arb on $x,y\in[0,12]$, with $t\ln t$ enclosed by its monotone pieces and $T\ge0$ by superadditivity of $t\ln t$, gives $0.25499\le\max T\le0.25551$; the tail $x+y>12$ contributes at most $12\ln2\,\mathrm e^{-12}<5\cdot10^{-5}$, since the bracket is at most $(x+y)\ln2$. The maximum is attained at $u=v=\mathrm e^{-1/2}\approx0.607$. For $g=M$ of Example~\ref{ex:matrix}, the explicit thresholds $u=(0.5876,0.8796,1-10^{-9})$ and $v=(10^{-9},0.4920,0.6812)$, evaluated in Arb, give the value $0.26253$, a rigorous lower bound for the coefficient of $M$ that exceeds the rigorous upper bound $0.25551$ for $\1\{M=4\}$. Thus the two tables have different first-order coefficients, although $\rho=1$ in both cases. For a general first-order term $T$ the first-order coefficient is still $\sup|\sum_{ij}\Delta_{ij}T(u_i,v_j)|$ over nondecreasing grids, which can be computed numerically. If $T=\sum_k\sigma_kA_k(u)B_k(v)$ with nonnegative unimodal $A_k,B_k$ normalised to maximum $1$, Lemma~\ref{lem:layer} gives the upper bound $\rho(g)\sum_k|\sigma_k|$, but in general no matching lower bound. For Gumbel such a decomposition exists as an integral: since $\frac{d^2}{dt^2}t\ln t=1/t=\int_0^\infty\mathrm e^{-wt}dw$,
\[
(x+y)\ln(x+y)-x\ln x-y\ln y=\int_0^x\!\!\int_0^y\frac{ds\,dt}{s+t}=\int_0^\infty(1-\mathrm e^{-wx})(1-\mathrm e^{-wy})\,\frac{dw}{w^2},
\]
so $T(u,v)=\int_0^\infty a_w(u)a_w(v)\,w^{-2}dw$ with the Huang--Kotz factors $a_w(u)=u(1-u^w)\ge0$, unimodal with $\max a_w=w(1+w)^{-1-1/w}$. Lemma~\ref{lem:layer} under the integral gives, for every $g$, a first-order coefficient at most $\rho(g)\int_0^\infty(1+w)^{-2-2/w}dw<0.31717\,\rho(g)$ (the integral lies in $[0.3171637,\,0.3171639]$ by interval arithmetic in Arb on $[10^{-7},10^9]$, the integrand being bounded by $1$ on $[0,10^{-7}]$ and by $(1+w)^{-2}$ beyond $10^9$), while the two-point construction of Proposition~\ref{prop:general} applied to $T\ge0$ gives at least $\rho(g)\max T\ge0.25499\,\rho(g)$. Thus $\rho$ determines the Gumbel coefficient up to the factor $0.31717/0.25499<1.25$, and the coefficients of $\1\{M=4\}$ and $M$ above lie in this range. The same argument applies to any Archimedean family whose $\chi''$ is completely monotone, $\chi''(x)=\int_{[0,\infty)}\mathrm e^{-tx}\mu(dt)$, provided $\int(1+t)^{-2-2/t}\mu(dt)<\infty$ (the integrand is $(\max a_t)^2/t^2$, read as $\mathrm e^{-2}$ at $t=0$): then $\chi(x+y)-\chi(x)-\chi(y)=\int_0^x\!\int_0^y\chi''(s+s')\,ds\,ds'=\int(1-\mathrm e^{-tx})(1-\mathrm e^{-ty})t^{-2}\mu(dt)$ by Fubini's theorem, $T=\int a_t(u)a_t(v)t^{-2}\mu(dt)$, and the first-order coefficient lies between $\rho(g)\max T$ and $\rho(g)\int(1+t)^{-2-2/t}\mu(dt)$. For the Joe family the expansion $\ln(1-\mathrm e^{-x})=-\sum_{m\ge1}\mathrm e^{-mx}/m$ gives $\chi(x)=-1+\sum_{m\ge1}\mathrm e^{-mx}/(m(m+1))$, so that $\chi''(x)=\sum_{m\ge1}\frac{m}{m+1}\mathrm e^{-mx}$ and $T(u,v)=\sum_{m\ge1}a_m(u)a_m(v)/(m(m+1))$, a countable mixture of the same Huang--Kotz products. The same argument gives first-order coefficients between $\rho(g)\max T\ge0.16235\,\rho(g)$ and $\rho(g)\sum_{m\ge1}m(1+m)^{-3-2/m}<0.19739\,\rho(g)$, so for Joe $\rho$ determines the coefficient up to a factor below $1.22$ (both constants certified in Arb: the series by a partial sum up to $m=2\cdot10^6$ with the tail bound $\sum_{m>N}(1+m)^{-2}<1/(N+1)$, $\max T$ from below by evaluation at $u=v=0.6909$; scripts \texttt{gumbel\_integral\_cert.py} and \texttt{joe\_constants\_cert.py} in the code archive).
\end{remark}

\section{EXAMPLES AND NUMERICAL EVIDENCE}\label{sec:numerics}

\begin{example}[An evaluation matrix]\label{ex:matrix}
The $4\times4$ aggregation matrix used as a worked example in \citet[p.~86, Fig.~1.34]{Alekseev2025} and its mixed differences are
\[
M=\begin{pmatrix}1&1&2&3\\1&2&2&3\\1&2&3&3\\2&3&3&4\end{pmatrix},\qquad
\Delta(M)=\begin{pmatrix}1&-1&0\\0&1&-1\\0&-1&1\end{pmatrix},
\]
with rows indexed by $X$. Here $N(M)=0$, $\rho(M)=1$ and $V(M)=6$. Compare this with the ``and'' matrix $m_{ij}=\min(i,j)$, for which $\Delta=I_3$ and $\rho=V=N=3$. Although $M$ has the larger total interaction $V$, its positive and negative interactions cancel over every rectangle. To first order, $M$ is therefore three times less sensitive to weak dependence than $\min$. The probability of the top grade, $g=\1\{M=4\}$, has a single nonzero mixed difference $\Delta_{33}=1$. By Proposition~\ref{prop:sign}, for $\lambda>0$ the independence computation never overestimates it, and the largest underestimation over all marginals is exactly $\arcsin(\lambda)/(2\pi)$, about $0.10$ at $\lambda=0.6$; for particular marginals it is smaller. For a designer choosing between $M$ and $\min$ as the aggregation rule, Table~\ref{tab:examples} quantifies the difference. At $\lambda=0.6$ an independence-based computation of the mean grade can be off by up to $0.31$ grade points for $\min$ (exactly, by Proposition~\ref{prop:sign}), but by only about $0.11$ for $M$. On a four-point scale this is the difference between an error of a third of a grade and one of a tenth.
\end{example}

Table~\ref{tab:examples} compares the bounds of Theorem~\ref{thm:gauss} with the numerically computed worst case. The worst case was obtained by maximising $|B|$ over the $2(r-1)$ thresholds with differential evolution (SciPy; seeds $0,1,2$; population size $30$, at most $400$ generations, tolerance $10^{-10}$; thresholds in $[-4,4]$) and taking the larger of the result and the lower bound. Any value found by optimisation is attained by an explicit choice of marginals, so, up to the floating-point error in evaluating the objective, it is a lower bound for $S_\lambda$; it need not be the global maximum. We write $\hat S_\lambda$ for the best value found.

\begin{table}[t]
\centering
\small
\caption{Worst-case bias of the independence computation under a Gaussian copula: lower bound $\frac{\arcsin\lambda}{2\pi}\rho$, best value found $\hat S_\lambda$, and the smallest of the upper bounds $\frac{\lambda}{2\pi}\rho+\min\{V R(\lambda),\rho G_2(\lambda)\}$ (Theorem~\ref{thm:gauss}) and $\frac{\arcsin\lambda}{2\pi}\max(P,Q)$ (Proposition~\ref{prop:sign}).}\label{tab:examples}
\begin{tabular}{@{}lcc ccc ccc@{}}
\toprule
& & & \multicolumn{3}{c}{$\lambda=0.2$} & \multicolumn{3}{c}{$\lambda=0.6$}\\
\cmidrule(lr){4-6}\cmidrule(l){7-9}
$g$ & $\rho$ & $V$ & lower & $\hat S_\lambda$ & upper & lower & $\hat S_\lambda$ & upper\\
\midrule
$\min(i,j)$ & 3 & 3 & 0.0961 & 0.0961 & 0.0961 & 0.3072 & 0.3072 & 0.3072\\
$M$ (Example~\ref{ex:matrix}) & 1 & 6 & 0.0320 & 0.0323 & 0.0374 & 0.1024 & 0.1108 & 0.1710\\
$\1\{M\ge3\}$ & 1 & 4 & 0.0320 & 0.0320 & 0.0366 & 0.1024 & 0.1024 & 0.1458\\
$\1\{M=4\}$ & 1 & 1 & 0.0320 & 0.0320 & 0.0320 & 0.1024 & 0.1024 & 0.1024\\
\bottomrule
\end{tabular}
\end{table}

\begin{table}[t]
\centering
\small
\caption{Sharpness of the remainder and of the lower bound. Largest per-cell remainder $\sup_{a,b}|c(a,b)-\lambda\varphi(a)\varphi(b)|$, computed on a grid, compared with $R(\lambda)$ in \eqref{eq:R} and with the Sz\'asz-based bound.}\label{tab:R}
\begin{tabular}{@{}lcccc@{}}
\toprule
$\lambda$ & $0.1$ & $0.3$ & $0.6$ & $0.9$\\
\midrule
largest per-cell remainder & 0.000294 & 0.002721 & 0.01230 & 0.03973\\
$R(\lambda)$ & 0.000294 & 0.002725 & 0.01258 & 0.04440\\
$(-\lambda-\ln(1-\lambda))/(2\pi)$ & 0.000853 & 0.009020 & 0.05034 & 0.2232\\
\bottomrule
\end{tabular}
\end{table}

To test the bounds systematically we enumerated all $21\,238$ monotone $4\times4$ matrices with entries in $\{1,\dots,4\}$, $m_{11}=1$ and $m_{44}=4$, and identified matrices whose $\Delta$ agree up to transposition and sign, which leaves $6535$ classes ($6534$ with $\rho>0$). All have $\max g-\min g=3$, and the design index $\rho(g)/(\max g-\min g)$ equals $0$, $1/3$, $2/3$ and $1$ for $1$, $190$, $2889$ and $3455$ classes, respectively; $M$ of Example~\ref{ex:matrix} is among the $190$ least sensitive. We then computed $\hat S_\lambda$ for $150$, $300$ and $150$ of these classes at $\lambda=0.2$, $0.6$ and $0.9$. Half of each sample was drawn at random (NumPy generator, seed $1$) among tables with $\rho>0$. The other half were the tables with the largest ratio $K(g)/\rho(g)$, where $K(g)=\frac{\arcsin\lambda}{2\pi}\max_U|\sum_{(i,j)\in U}\Delta_{ij}(g)|$ and $U$ runs over the shapes of superlevel sets of $(i,j)\mapsto c(a_i,b_j)$ observed for $6\times10^4$ random threshold vectors at each of three scales (seed $0$). If all possible shapes were included, the layer-cake argument applied to $c$ would give $S_\lambda\le K$. A large ratio $K/\rho$ therefore flags tables for which two-point marginals may be suboptimal. Each table was optimised with two differential-evolution seeds and the same settings as above. The upper bound in \eqref{eq:main} was never violated. Its ratio to $\hat S_\lambda$ had median $1.14$, $1.36$, $1.69$ and maximum $1.33$, $2.04$, $3.05$, respectively. The slack comes from summing per-cell remainders with weight $V(g)$, which ignores cancellations. The smaller of the two upper bounds, $\frac{\lambda}{2\pi}\rho+\min\{VR(\lambda),\rho G_2(\lambda)\}$, was never violated either; its ratio to $\hat S_\lambda$ had maximum $1.17$, $1.77$ and $3.05$, and the tetrachoric bound was the smaller one for $30\%$, $3\%$ and $0\%$ of the sampled tables ($V/\rho$ ranged from $1$ to $9$). Thus at moderate dependence ($\lambda\approx0.6$) the upper bound of Theorem~\ref{thm:gauss} can exceed the true worst case by a factor of up to about $1.8$, and by up to $3$ at $\lambda=0.9$. It is best read as an asymptotic statement with an explicit error term. For a numerical value at a given $\lambda$, the worst case itself can be computed by optimisation over the $(r-1)+(c-1)$ thresholds. Because optimisation yields lower bounds for $S_\lambda$, the reported excesses over the lower bound of Theorem~\ref{thm:gauss} (up to $30\%$ at $\lambda=0.9$) are attained by explicit marginals, up to floating-point evaluation error. The comparison with the upper bound is only a consistency check.

To check the optimisation independently, for $\Delta$ with $+1$ on the diagonal and $-1$ on the superdiagonal we compared differential evolution (five seeds, population $40$, polishing) with $300$ random starts of L-BFGS-B. At $\lambda=0.6$ the best values were $0.1140051251$ and $0.1140051251$, and at $\lambda=0.9$ they were $0.2325117233$ and $0.2325117233$; they agree to $10^{-11}$, and so do the ten best L-BFGS-B runs. For this table the worst marginals found have full support on the four grades: at $\lambda=0.6$, $p_X=(0.157,0.261,0.213,0.369)$, with $p_Y$ the same vector reversed; at $\lambda=0.9$, $p_X=(0.182,0.262,0.271,0.285)$, likewise reversed. This is a numerical observation for one table, not a general property, and agreement of two optimisers does not prove global optimality.

Two-point marginals are asymptotically optimal but need not be optimal at finite dependence. Over these samples, the largest ratio $\hat S_\lambda(g)\big/\bigl(\frac{\arcsin\lambda}{2\pi}\rho(g)\bigr)$ was $1.015$, $1.113$ and $1.305$ at $\lambda=0.2$, $0.6$ and $0.9$. At $\lambda=0.6$ some tables admit marginals with more support points that exceed the lower bound by $11\%$; an example is $\Delta$ with $+1$ on the diagonal and $-1$ on the superdiagonal. At $\lambda=0.9$ the excess reaches $30\%$. The reason is that $c(a,b)$ is no longer close to a product: its superlevel sets on the threshold grid take non-rectangular shapes, so the layer-cake argument of Lemma~\ref{lem:layer} no longer applies. Hoeffding's reduction to few-point laws is thus exact for the first-order term and can fail beyond it.

The worst case is a guarantee, not a typical value. For uniform marginals on four grades the biases of $\min(i,j)$, $M$ and $\1\{M=4\}$ are $0.0744$, $-0.0019$ and $0.0212$ at $\lambda=0.2$, that is $77\%$, $6\%$ and $66\%$ of the best values found, and $0.2464$, $0.0106$ and $0.0720$ at $\lambda=0.6$.

\begin{example}[Staging of colon cancer]\label{ex:colon}
In the TNM system the stage of a tumour is read from a table of two ordinal grades, the depth of invasion $T$ and the nodal category $N$. For node-positive (stage~III) colon cancer the 7th edition of the AJCC manual \citep{Edge2010,Jessup2011} assigns the substages IIIA, IIIB and IIIC, coded $1,2,3$. On the grades $T\in\{\mathrm{T1},\mathrm{T2},\mathrm{T3},\mathrm{T4b}\}$ (rows) and $N\in\{\mathrm{N1},\mathrm{N2a},\mathrm{N2b}\}$ ($1$--$3$, $4$--$6$ and at least $7$ positive nodes; columns) the table and its mixed differences are
\[
g=\begin{pmatrix}1&1&2\\1&2&2\\2&2&3\\3&3&3\end{pmatrix},\qquad
\Delta(g)=\begin{pmatrix}1&-1\\-1&1\\0&-1\end{pmatrix},
\]
so $\rho(g)=1$ and $V(g)=5$. For the share of substage IIIC, $g_C=\1\{g=3\}$, the matrix $\Delta(g_C)$ has a single entry $+1$ and a single entry $-1$, so $\rho=P=Q=1$. Registries and trial reports often tabulate $T$ and $N$ separately, and the substage mix is then obtained from the two margins as if $T$ and $N$ were independent; the supremum over marginals is exactly what is needed when the same table is applied to populations with different case mix.

We use the data set \texttt{colon} of the R package \texttt{survival} \citep{Therneau2026}, from the trial of adjuvant chemotherapy of \citet{Laurie1989} and \citet{Moertel1990}: $909$ patients with at least one positive node and a recorded node count, the extent of local spread (submucosa, muscle, serosa, adjacent structures) being read as T1, T2, T3, T4b. The counts are $(18,2,0)$, $(78,19,9)$, $(470,145,127)$ and $(27,7,7)$ for T1, T2, T3 and T4b across N1, N2a, N2b. The polychoric correlation is $\hat\lambda=0.152$ (bootstrap $95\%$ interval $[0.050,\,0.248]$, $400$ resamples of the patients with replacement, with $\lambda$ and the errors reported below recomputed in each, the margins being fixed at each resample's own cumulative proportions), and the Gaussian copula fits the table ($G^2=8.1$ on $5$ degrees of freedom, $p=0.15$).

Since the Gaussian copula with $\lambda>0$ is positively quadrant dependent, Propositions~\ref{prop:general} and~\ref{prop:sign} give $S_\lambda(g_C)=\varepsilon(\lambda)=\arcsin(\lambda)/(2\pi)$ exactly. At $\hat\lambda$ this is $0.0243$: whatever the distributions of $T$ and $N$, the independence computation misstates the share of IIIC by at most $2.4$ percentage points ($4.0$ at the upper end of the interval for $\lambda$). In these data the share of IIIC, as reconstructed from the recorded categories, is $18.5\%$ and the independence computation gives $17.4\%$, an underestimate of $1.13$ points (bootstrap $95\%$ interval $[0.32,\,1.91]$), below the bound. The identity $S_C(g_C)=\varepsilon(C)$ does not rest on the Gaussian model: by Proposition~\ref{prop:general} it holds for every positively quadrant dependent copula $C$, so the guarantee for the share of IIIC depends on the latent copula only through $\varepsilon(C)=\sup|C(u,v)-uv|$, and the Gaussian fit serves only to estimate it: the value $0.0243$ is the estimate of $\varepsilon(C)$ under the fitted Gaussian model, and under another model of the latent dependence $\varepsilon(C)$ would have to be estimated within that model; for an arbitrary copula the same proposition still gives $\varepsilon(C)\le S_C(g_C)\le2\varepsilon(C)$, and Theorem~\ref{thm:tv} gives $S_C(g_C)\le\tau(C)/2$. The estimate does depend on the model. Fitting one-parameter families by maximum likelihood, with the margins fixed at the observed cumulative proportions as for $\hat\lambda$, the Clayton ($\hat\theta=0.254$, $G^2=5.3$, $p=0.38$) and Frank ($\hat\theta=0.973$, $G^2=8.7$, $p=0.12$) copulas, which also fit the table, give $\varepsilon=0.031$ and $0.030$, that is, guarantees of $3.1$ and $3.0$ points, whereas the Gumbel copula ($\hat\theta=1.052$), which fits worse ($G^2=13.2$, $p=0.02$), gives $1.3$. These are point estimates; refitting all four families and recomputing $\varepsilon(\hat C)$ in each of the $400$ bootstrap resamples used for $\hat\lambda$ gives $95\%$ percentile intervals $[1.3,\,4.9]$, $[0.8,\,5.1]$ and $[0,\,3.0]$ points for the guarantee, and $[0.8,\,4.0]$ for the Gaussian copula. The observed point error is below all four point guarantees; separate intervals do not show how reliably it lies below them. This is a consistency check on one cohort, not a validation of the bounds for another population. For the substage index $g$, Theorem~\ref{thm:gauss} gives $0.0243\le S_{\hat\lambda}(g)\le\rho(g)G(\hat\lambda)=0.0273$, whereas the bound $V(g)\varepsilon$ would give $0.122$; the observed error is $-0.002$, with a bootstrap interval $[-0.011,\,0.007]$ that contains zero. The conclusions do not change if the $12$ records whose node count disagrees with the recorded indicator of more than four positive nodes are dropped ($\hat\lambda=0.149$; error for IIIC $1.13$ points), or if the cell (serosa, N2a), which may contain T4a tumours that the 7th edition classifies as IIIC, is coded IIIC ($\rho$ unchanged; error $1.55$ points). The computations are in \texttt{colon\_example.py} in the code archive.
\end{example}

\section{DISCUSSION}\label{sec:discussion}

\textbf{A dependence-free design index.} Theorems~\ref{thm:gauss} and~\ref{thm:families} show that, under weak dependence, the worst-case bias splits into a model factor (the strength $|\theta-\theta_0|$ and the constant $\kappa$ of the family) and a table factor $\rho(g)$. The table factor is the same for all families with product-form first-order departure. When aggregation tables are designed, $\rho(M)$ and $\rho(\1\{M\ge z\})$ for the relevant levels $z$ can therefore serve as robustness criteria, to be traded off against other requirements on $M$: $\rho(g)$ when $g$ is a numerical score whose mean is meaningful, and $\rho_{\mathrm{TV}}(g)=\max_A\rho(\1\{g\in A\})$, the first-order index of the total variation distance (see below), when only the ordinal grade matters. The index is cheap to compute: Kadane's algorithm applied to all pairs of row bounds gives $\rho$ in $O(r^2c)$ time \citep{Bentley1984}. The supremum over all marginals makes $\rho$ a guarantee rather than a typical value; restricting the marginals, for instance to $\Prob(X=i)\ge p_0$, can only lower the worst case, and which index then replaces $\rho$ is open. Proposition~\ref{prop:sign} also shows that a table with nonzero corner contrast $N$ cannot have worst-case bias below $|N|\arcsin(|\lambda|)/(2\pi)$, whatever else is done. For the Gaussian copula the role of $\rho$ extends beyond first order, since $S_\lambda(g)\le\rho(g)G(|\lambda|)$; the Vitali variation enters only the other remainder bound. Since all bounds increase in $|\lambda|$, an upper confidence limit for $|\lambda|$, for instance from the polychoric maximum-likelihood estimator and its asymptotic standard error \citep{Olsson1979}, gives an upper bound for $S_\lambda(g)$ that holds with at least the same confidence. This is useful when a table is designed before data are collected or transferred to populations without joint data; when joint data on $(X,Y)$ are available, the bias can be estimated directly.

\textbf{Distribution of the merged grade.} The theorems bound the bias of $\E g(X,Y)$, but they also control the whole distribution of the merged grade. Since the supremum over marginals commutes with a maximum, the worst-case Kolmogorov distance between the laws of $g(X,Y)$ under $C$ and under independence is $\max_zS_C(\1\{g\ge z\})$, and the worst-case total variation distance is $\max_AS_C(\1\{g\in A\})$, the maximum being over sets $A$ of values of $g$. Theorems~\ref{thm:gauss} and~\ref{thm:families} apply to each of the finitely many terms, so for families with a product-form first-order term the worst-case total variation distance is $|\theta-\theta_0|\kappa\rho_{\mathrm{TV}}(g)+O((\theta-\theta_0)^2)$ and the worst-case Kolmogorov distance is $|\theta-\theta_0|\kappa\max_z\rho(\1\{g\ge z\})+O((\theta-\theta_0)^2)$; for the Gaussian copula the first-order coefficients are $\max_z\rho(\1\{g\ge z\})/(2\pi)$ and $\max_A\rho(\1\{g\in A\})/(2\pi)$. For a monotone table $g$, each set $\{g\ge z\}$ is an upper set, and the $g$-volume form of $\rho$ shows that $\rho(\1\{g\ge z\})\in\{0,1\}$; since a $0$--$1$ table is additive only if it depends on one grade alone, and nested upper sets cannot switch between the two grades, all non-constant monotone tables that depend on both grades are equally sensitive in the Kolmogorov distance to first order, as are $M$ of Example~\ref{ex:matrix} and $\min(i,j)$ on four grades. But $\rho(\1\{M=3\})=2$, whereas $\rho(\1\{\min(i,j)\in A\})\le1$ for every $A$, so in total variation $M$ is twice as sensitive as $\min$ to first order, although for the mean grade it is three times less sensitive.

\textbf{Relation to ordinal latent-normal models.} If $\Delta(g)=\alpha\beta^\top$ with $\alpha,\beta\ge0$, the case of the covariance of two increasing score functions, the optimal rectangle is the whole grid and the optimal marginals split each variable at the median. This agrees with the elementary fact that the orthant covariance $\Phi_2(a,b;\lambda)-\Phi(a)\Phi(b)=\int_0^\lambda\varphi_2(a,b;s)\,ds$ of a dichotomised bivariate normal pair is largest at the median split $a=b=0$, where it equals $\arcsin(\lambda)/(2\pi)$ \citep{Sheppard1899}. The new element is the mixed-sign case, where the optimal rectangle may be a proper contiguous block.

\textbf{Related extremal problems.} By \eqref{eq:identity}, $B=\iint(C-\Pi)\,dm$ for the signed measure $m=\sum_{i,j}\Delta_{ij}\,\delta_{(u_i,v_j)}$ on $[0,1]^2$. For $\Delta=\alpha\beta^\top$ with $\alpha,\beta\ge0$ the measure factorises, $m=\mu\otimes\nu$ with $\mu=\sum_i\alpha_i\delta_{u_i}$ and $\nu=\sum_j\beta_j\delta_{v_j}$, and $B$ is the covariance of two nondecreasing score functions. This is the setting of the Hoeffding-type bounds of \citet{Cuadras2002} and \citet{Beare2009}. In the mixed-sign case the extremum is governed by $\rho$, an interval-restricted analogue of the cut norm $\max_{I,J}|\sum_{I\times J}\Delta_{ij}|$ over arbitrary index sets \citep{FriezeKannan1999,AlonNaor2006}. The interval restriction is what makes $\rho$ computable exactly in polynomial time, whereas computing the cut norm is MAX SNP hard, so it has no polynomial-time approximation scheme unless P${}={}$NP, although it can be approximated in polynomial time within a factor $0.56$ \citep{AlonNaor2006}.

\textbf{Second-order term.} The $O(\lambda^2)$ term in Theorem~\ref{thm:gauss} is probably not sharp in order. The second-order part of $c(a,b)$ is $\frac{\lambda^2}2\,a\varphi(a)\,b\varphi(b)$, which vanishes at the median split $a=b=0$ used in the lower bound. To test the order we computed $\hat S_\lambda(g)$ for $g=M$ of Example~\ref{ex:matrix} and for $\Delta$ with $+1$ on the diagonal and $-1$ on the superdiagonal at $\lambda=0.02,0.04,0.08,0.16,0.32$. The orthant covariances were evaluated by $40$-point Gauss--Legendre quadrature of $\int_0^\lambda\varphi_2\,ds$, and the thresholds were optimised by differential evolution followed by Nelder--Mead. On a log--log scale the local slopes of $\hat S_\lambda(g)-|\lambda|\rho(g)/(2\pi)$ were $2.98$--$3.00$ and $2.92$--$2.95$ for the two tables, and those of the excess $\hat S_\lambda(g)-\frac{\arcsin|\lambda|}{2\pi}\rho(g)$ over the lower bound were $2.96$--$2.99$ and $2.88$--$2.93$. For the two tables examined, the numerical results are thus consistent with a cubic departure from two-point optimality. \emph{Conjecture:} $S_\lambda(g)=|\lambda|\rho(g)/(2\pi)+O(|\lambda|^3)$ whenever the rectangle attaining $\rho(g)$ in \eqref{eq:rho} is unique. Among the $6534$ classes of $4\times4$ tables with $\rho>0$ this holds for $1138$, about $17\%$; ties often come from rows or columns of zeros in $\Delta$ next to the optimal rectangle. This is supported by the numerical slopes above but not proved. The conjecture holds when $\Delta(g)$ has constant sign, by Proposition~\ref{prop:sign}, in particular for all $2\times2$ tables. In general, by the proof of Theorem~\ref{thm:gauss}, the terms with $n\ge3$ of the tetrachoric series contribute at most $\rho(g)\sum_{n\ge3}|\lambda|^nT_n^2/(4\,n!)=O(\rho(g)|\lambda|^3)$ uniformly in the marginals, so the conjecture concerns only the first two terms; then the maximiser of $|T_1|$ below is the median split, where the second-order term vanishes. If several nested rectangles attain $\rho(g)$, the maximisers of $|T_1|$ include thresholds away from $0$, and we cannot exclude a nonzero $\lambda^2$ term for such degenerate tables. The remainder bounds $V(g)R(|\lambda|)$ and $\rho(g)G_2(|\lambda|)$, both of order $\lambda^2$, are therefore a limitation of the method, not a proven property of general tables. This conjecture is consistent with the lower bound, whose excess $(\arcsin|\lambda|-|\lambda|)\rho/(2\pi)$ is of order $|\lambda|^3$. A proof needs a finer upper bound than \eqref{eq:main}, because the per-cell second-order term is genuinely of order $\lambda^2$ away from the median. Writing $B=\lambda T_1+\frac{\lambda^2}{2}T_2+O(\lambda^3)$ as functions of the thresholds, it would suffice that $T_2$ vanishes on the set of maximisers of $|T_1|$ and that $|T_2|$ grows at most linearly, while $\rho/(2\pi)-|T_1|$ grows at least quadratically, with the distance from that set.

\section*{DECLARATIONS}
\small
\textbf{Code availability.} The code for all computations and computer-assisted bounds is available at Zenodo, \url{https://doi.org/10.5281/zenodo.24326496}.

\textbf{Funding.} This research received no specific grant from any funding agency in the public, commercial, or not-for-profit sectors.

\textbf{Competing interests.} The author declares no competing interests.

\textbf{Generative AI.} During the preparation of this work the author used Claude (Anthropic) in order to assist with the literature search, with the typesetting and presentation of mathematical derivations, and with the implementation of the computer code. After using this tool, the author reviewed and edited the content as needed and takes full responsibility for the content of the publication.

\makecontacts

\end{document}